\documentclass[11pt,reqno]{amsart}
\usepackage[T1]{fontenc}
\usepackage{lmodern}
\usepackage{amsmath,amssymb,amsthm,mathtools}
\usepackage[letterpaper,margin=1.15in]{geometry}
\usepackage{microtype}
\usepackage[hidelinks]{hyperref}
\hypersetup{pdftitle={A counterexample to Rubio de Francia's weighted Littlewood--Paley conjecture},pdfauthor={Veli Shakhmurov and Rishad Shahmurov},pdfsubject={Weighted harmonic analysis},pdfkeywords={Rubio de Francia, square functions, A1 weights, Fourier projections}}
\allowdisplaybreaks[2]
\numberwithin{equation}{section}
\newtheorem{theorem}{Theorem}[section]
\newtheorem{lemma}[theorem]{Lemma}
\newtheorem{corollary}[theorem]{Corollary}
\theoremstyle{remark}
\newtheorem{remark}[theorem]{Remark}
\newcommand{\R}{\mathbb R}
\newcommand{\T}{\mathbb T}
\newcommand{\Z}{\mathbb Z}
\newcommand{\1}{\mathbf 1}
\newcommand{\dd}{\,\mathrm d}
\newcommand{\supp}{\operatorname{supp}}
\newcommand{\dist}{\operatorname{dist}}
\newcommand{\norm}[1]{\lVert #1\rVert}
\newcommand{\abs}[1]{\lvert #1\rvert}
\DeclareMathOperator{\esssup}{ess\,sup}
\DeclareMathOperator{\essinf}{ess\,inf}
\title[The weighted Rubio de Francia conjecture]{A counterexample to Rubio de Francia's weighted Littlewood--Paley conjecture}
\author[V. Shakhmurov]{Veli Shakhmurov}
\address{Department of Industrial Engineering, Antalya Bilim University, Antalya 07190, Turkey; Azerbaijan State Economic University, Baku AZ1001, Azerbaijan}
\email{veli.sahmurov@antalya.edu.tr}
\email{veli.sahmurov@gmail.com}
\author[R. Shahmurov]{Rishad Shahmurov}
\address{Cellular Products Research and Development, 595 East Crossville Road, Roswell, GA 30075, USA}
\email{rshahmurov@crimson.ua.edu}
\thanks{Corresponding author: Rishad Shahmurov.}
\date{November 27, 2026}
\begin{document}
\begin{abstract}
Rubio de Francia's square-function theorem extends Littlewood--Paley theory to arbitrary disjoint frequency intervals. We disprove his 1985 conjecture for the strong $L^2(w)$ estimate with $w\in A_1$. For any prescribed finite dilation separation, the frequency family can be fixed, dyadic, and contained in a bounded window. The estimate fails for every positive power of a single even weight in $A_1$, whose small powers have characteristic tending to one. The proof amplifies a small spectral component concentrated away from the input support into almost complete transfer of square-function energy to a set of arbitrarily small maximal density. The obstruction persists for smooth square functions and finite maximal-function bumps.
\end{abstract}
\maketitle

\section{Introduction and main results}
We use the Fourier normalization
\[
\widehat f(\xi)=\int_{\R} f(x)e^{-2\pi i x\xi}\dd x,
\qquad S_\omega f=\mathcal F^{-1}(\1_\omega\widehat f).
\]
For a finite or countable collection $\Omega$ of pairwise disjoint bounded intervals of $\R$, set
\[
T_\Omega f(x)=\left(\sum_{\omega\in\Omega}\abs{S_\omega f(x)}^2\right)^{1/2}.
\]
In particular, Plancherel's theorem gives
\begin{equation}\label{eq:plancherel}
\norm{T_\Omega f}_2^2=\sum_{\omega\in\Omega}\norm{S_\omega f}_2^2\le\norm f_2^2.
\end{equation}
Let $M$ be the uncentered Hardy--Littlewood maximal operator on $\R$. A weight is a locally integrable function $w>0$ almost everywhere, and
\[
[w]_{A_1}=\esssup_{x\in\R}\frac{Mw(x)}{w(x)}.
\]
\subsection{The endpoint in the original question}
Classical Littlewood--Paley theory compares a function's $L^p$ norm with the square function of its dyadic frequency pieces. For arbitrary disjoint intervals, Rubio de Francia \cite{RF85} established the upper bound for $2\le p<\infty$; the range $1<p<2$ fails. His theorem unified lacunary and congruent frequency decompositions. The introduction of \cite[p.~2]{RF85} credits C\'ordoba with bringing the general question to his attention and notes the application of the earlier congruent-interval estimate to Bochner--Riesz multipliers.

Rubio de Francia also proved the weighted bound for $p>2$ and $w\in A_{p/2}$, with a sharp weight-class index \cite[Theorem~6.1]{RF85}. He posed the endpoint for two reasons: it is the limiting case of that theorem, and it holds for the two extremal frequency geometries (1.4) and (1.5) of his paper. The proposed estimate \cite[Section~6, (6.4), p.~10]{RF85} is
\begin{equation}\label{eq:conjecture}
\int_{\R}T_\Omega f(x)^2w(x)\dd x
\le C_w\int_{\R}\abs{f(x)}^2w(x)\dd x,
\qquad w\in A_1,
\end{equation}
for arbitrary disjoint interval families. The question is also discussed in \cite[Section~8.2]{Duo01}.

The endpoint would provide a form of Fourier orthogonality robust under spatial reweighting. Indeed, for a finite family and independent random signs $\epsilon_\omega$,
\[
\mathbb E_\epsilon\int_{\R}\Big|\sum_{\omega\in\Omega}\epsilon_\omega S_\omega f\Big|^2w
=\int_{\R}T_\Omega f^2w.
\]
Thus \eqref{eq:conjecture} would control the mean squared norm of randomized interval multipliers. The unweighted square-function theorem already has important consequences: Coifman, Rubio de Francia and Semmes \cite{CRS88} used it for multipliers of finite $q$-variation on dyadic frequency intervals ($2\le q<\infty$, $1<p<\infty$, $|1/p-1/2|<1/q$). Lacey \cite{Lacey07} surveys the time-frequency proofs, higher-dimensional extensions and multiplier applications. The weighted endpoint would give additional control of frequency-resolved energy in physical space.

\subsection{Results}
The following result gives a negative answer to \eqref{eq:conjecture}.
Write $\mathcal D_{\rm dyad}=\{[m2^k,(m+1)2^k):m,k\in\Z\}$ for the standard dyadic grid in frequency. For $\Lambda\ge1$, $\Lambda\omega$ denotes the concentric interval of length $\Lambda|\omega|$. We will require the dilated intervals to have disjoint closures.
\begin{theorem}[A fixed family and a fixed weight direction]\label{thm:weighted}
For every prescribed $\Lambda\ge1$ there are an even weight $0<W\le1$ in $A_1(\R)$ and a single countable family $\Omega_*\subset\mathcal D_{\rm dyad}$ contained in $(0,1)$, with pairwise disjoint closed dilations $\Lambda\omega$, such that, for every $t>0$,
\begin{equation}\label{eq:unbounded}
\sup_{0\ne f\in C_c^\infty(\R)}
\frac{\norm{T_{\Omega_*}f}_{L^2(W^t)}}{\norm f_{L^2(W^t)}}=\infty.
\end{equation}
For $0<t\le1$, $[W^t]_{A_1}\le[W]_{A_1}^t$. In particular, for every $\varepsilon>0$ a counterexample weight with characteristic at most $1+\varepsilon$ is obtained from this one weight and family. The endpoint fails even if its constant is allowed to depend on both the weight and the family.
\end{theorem}

Our main tool is the following concentration result. Its geometric hypothesis is stronger than merely requiring $E$ to be small in measure: it imposes a density restriction in every interval containing an input point.
\begin{theorem}[Thin-set concentration]\label{thm:thin}
Given $0<\delta,\kappa<1$ and $\Lambda\ge1$, there exist bounded measurable sets $E,F\subset\R$, a nonzero $f\in C_c^\infty(\R)$ with $\supp f\subset F$, and a finite family $\Omega\subset\mathcal D_{\rm dyad}$ with pairwise disjoint closed dilations $\Lambda\omega$, such that
\begin{equation}\label{eq:thin}
M\1_E(y)\le\delta\quad (y\in F),
\qquad \int_E T_\Omega f(x)^2\dd x>\kappa\norm f_2^2.
\end{equation}
For given $\delta,\Lambda$, call a configuration satisfying these boundedness, support, frequency and maximal-density requirements admissible. The supremum of $\norm{\1_E T_\Omega f}_2^2/\norm f_2^2$ over all such configurations is $1$.
\end{theorem}

\noindent\emph{Coherent and incoherent energy.} In any example of Theorem~\ref{thm:thin}, let $P_\Omega=\sum_{\omega\in\Omega}S_\omega$. This is an orthogonal projection on unweighted $L^2$. Since $f=0$ almost everywhere on $E$,
\[
\norm{\1_E P_\Omega f}_2^2\le\norm{P_\Omega f-f}_2^2
=\norm f_2^2-\norm{P_\Omega f}_2^2<(1-\kappa)\norm f_2^2.
\]
Here we used $\norm{P_\Omega f}_2^2=\norm{T_\Omega f}_2^2\ge\norm{\1_E T_\Omega f}_2^2$. Thus the squared moduli of the individual projections can account for nearly all the input energy on $E$, while their coherent sum has arbitrarily little energy there. The difference lies in the cross terms, not in a violation of Plancherel.

\subsection{Partial results and the difficulty of the endpoint}
Congruent intervals were treated in \cite{RF83}; lacunary decompositions are governed by classical weighted Littlewood--Paley theory. The power weights $|x|^{-\alpha}$, $0<\alpha<1$, were treated by Rubio de Francia \cite{RF89} and by an alternative argument of Carbery \cite{Carbery91}. Garg, Roncal and Shrivastava \cite[Section~4.2]{GRS21} also identify a mixed case: congruent subdivisions within lacunary parent intervals. These positive cases involve either a structured frequency geometry or a spatial weight with controlled shape.

Kr\'ol \cite[Theorem~B]{Krol14} proved the weighted weak endpoint and strong estimates for larger sequence exponents (see Corollary~\ref{cor:q}); Garg, Roncal and Shrivastava \cite{GRS21} obtained quantitative bounds by sparse and time-frequency methods. Di Plinio, Fl\'orez-Amatriain, Parissis and Roncal \cite{DFPR25} improved the weak bound and proved the strong endpoint for the Walsh model and even radially decreasing weights. Their Section~1.3, p.~449, records the general Fourier endpoint as open.

Berkson's even decreasing-weight result \cite{Berkson23} is compatible with our construction; his earlier affirmative claims \cite{Berkson22} are examined in Appendix~\ref{app:literature}. The abstract of his later paper \cite{Berkson24} reports a related estimate involving a rearranged input weight; no result from that paper is used here.

The original sharp-function estimate, controlled by $(M|f|^2)^{1/2}$, is effective above $p=2$ but does not yield strong control at the endpoint. The distinct conjecture with sequence exponent $p'$ below two \cite[Conjecture~7.2, p.~10]{RF85} was disproved by Cowling and Tao \cite{CT01}.

\subsection{How the construction arises}
The upper bound in Theorem~\ref{thm:thin} is just Plancherel. The lower bound reveals what the weighted problem is asking for. If $E$ is thin as seen from $F$, then the standard maximal weight $v=(M\1_E)^\alpha$, $0<\alpha<1$, has uniformly bounded $A_1$ characteristic, equals one almost everywhere on $E$, and is at most $\delta^\alpha$ on $F$. An affirmative uniform endpoint would therefore force
\[
\int_E T_\Omega f^2\lesssim_\alpha\delta^\alpha\norm f_2^2
\qquad(\supp f\subset F).
\]
This suggests trying to move Fourier energy into thin receiving regions. The issue is not the absolute size of a single projection, but the fraction of \emph{its own} energy lying outside the input support.

A single projection suggests the first step. Construct a normalized function $g$ supported in a small interval $H$ whose low-frequency part $g_0$ has almost all of its own energy in a disjoint interval $V$. Finite moment interpolation allows us to prescribe this low-frequency part. The cost is that $p=\norm{g_0}_2^2$ may be tiny.

To amplify $p$, multiply by a rapidly oscillating copy of $g$. In the new low-frequency block retain the whole preceding product; in each new high-frequency block repeat the preceding selection. If $\gamma=\int_V|g_0|^2/p$, the two steps yield $\gamma p+(1-p)\gamma p$. Repetition gives the ideal receiving energy
\begin{equation}\label{eq:amplify}
\gamma p\sum_{j=1}^N(1-p)^{N-j}=\gamma[1-(1-p)^N].
\end{equation}
Choose $\gamma$ close to one, then $N$ large. Every frequency is selected at most once, so this amplification respects Plancherel.

The spatial recursion places receivers at the first failed support condition. Disjoint receiving cells prevent density from accumulating. The construction thus rests on a general principle: a spectral component may carry very little of the total energy of a function while exhibiting almost complete \emph{relative} concentration away from its support. Recursive frequency selection amplifies this small component without violating Plancherel orthogonality.

For the guarded construction developed here, the optimal one-factor energy is exponentially small in the inverse maximal density: Appendix~\ref{app:cost} proves $p_* = \exp(-\Theta(1/\delta))$ for fixed input interval, concentration accuracy and dilation separation. For a fixed target fraction below the one-factor concentration guarantee, the least depth prescribed by the geometric-series amplification has order $N=\exp(\Theta(1/\delta))$. This quantitative feature helps explain why the obstruction may remain invisible in moderate-depth constructions or numerical experiments. It concerns the present amplification scheme, not a universal lower bound for arbitrary counterexamples. The principle also suggests a method for investigating other endpoint problems in weighted Fourier analysis; Section~\ref{sec:scope} gives concrete extensions.

Sections~\ref{sec:moments}--\ref{sec:energy} make the construction precise; Section~\ref{sec:weight} fixes the weight and the family. The contrast formula and the limiting case of Kr\'ol's sequence-exponent estimates follow in Section~\ref{sec:thresholds}. Section~\ref{sec:elliptic} transfers the obstruction to arbitrary spectral intervals of elliptic operators by a local rescaling argument.

\section{A small spectral component concentrated off the input support}\label{sec:moments}
Fix the desired dilation factor $\Lambda\ge1$ throughout the construction.
We identify $\T=\R/\Z$ with $[0,1)$, equipped with normalized Lebesgue measure, and write $\widehat g(n)=\int_0^1 g(t)e^{-2\pi int}\dd t$. For a real interval $\omega$, let
\[
S_\omega^{\T}g(t)=\sum_{n\in\omega\cap\Z}\widehat g(n)e^{2\pi int}.
\]
The notation $S_\omega^{\R}$ will be used when it is necessary to distinguish Fourier integrals from Fourier series.

\begin{lemma}[Prescribed moments]\label{lem:moment}
Let $H,V\subset\T$ be nonempty open intervals with disjoint closures. For any $0<\eta<1$ there exist a power of two $L\ge4$, an integer $G=\lceil(\Lambda+2)L\rceil$, and $g\in C_c^\infty(H)$ with $\norm g_{L^2(\T)}=1$ such that, with
\[
g_0=S_{[-L,L]}^{\T}g,\qquad g_1=g-g_0,\qquad p=\norm{g_0}_{L^2(\T)}^2,
\]
one has
\begin{equation}\label{eq:moment}
0<p<1,\qquad \int_V\abs{g_0(t)}^2\dd t>(1-\eta)p,\qquad \widehat g(n)=0\quad (L/2-1<\abs n\le G).
\end{equation}
\end{lemma}
\begin{proof}
We first prescribe the desired low-frequency output, rather than try to guess the input. Let $v$ be the midpoint of $V$, and take the normalized Dirichlet polynomial
\[
b(t)=(2l+1)^{-1/2}\sum_{|n|\le l}e^{2\pi in(t-v)}.
\]
It has $\norm b_2=1$. For $s=t-v$ represented in $[-1/2,1/2]$ we have $|b(t)|^2\le ((2l+1)4s^2)^{-1}$ when $s\ne0$. Consequently
\[
\int_{\T\setminus V}|b|^2\le \frac{1}{(2l+1)|V|}.
\]
Choosing $l\ge1$ with $(2l+1)|V|>\eta^{-1}$ ensures $\int_V|b|^2>1-\eta$. Write $b(t)=\sum_{|n|\le l}b_ne^{2\pi int}$. Let $L$ be the smallest power of two at least $2l+2$, set $G=\lceil(\Lambda+2)L\rceil$, and extend the coefficient vector by zero for $l<|n|\le G$. This empty frequency region leaves room for smooth cutoffs and for separated dilations. Choose $t_0\in H$ and a nonnegative, nonzero bump $\phi\in C_c^\infty(H)$ supported so close to $t_0$ that $2\pi G|t-t_0|<\pi/3$ on its support. For $|n|\le G$ the real part of $e^{2\pi int_0}\widehat\phi(n)$ is positive, so these coefficients are nonzero. Ordinary polynomial interpolation gives a polynomial $Q$ of degree at most $2G$ satisfying
\[
Q(2\pi in)=b_n/\widehat\phi(n),\qquad -G\le n\le G.
\]
Set $u=Q(\mathrm d/\mathrm dt)\phi$. Differentiation preserves support, and $\widehat u(n)=Q(2\pi in)\widehat\phi(n)=b_n$ for $|n|\le G$. Thus $g=u/\norm u_2$ gives $g_0=b/\norm u_2$ and $p=\norm u_2^{-2}>0$. If $p=1$, orthogonality implies $g=g_0$; a nonzero trigonometric polynomial cannot be supported in $H$. Thus $p<1$, and the concentration of $b$ and the extra zero moments prove \eqref{eq:moment}.
\end{proof}

No quantitative bound on the inverse moment map is needed. The bandwidth has the right order: for any nonzero polynomial $b$ of degree $l$, Cauchy--Schwarz gives $\norm b_\infty^2\le(2l+1)\norm b_2^2$, so concentration of a fraction $1-\eta$ in $V$ requires $(2l+1)|V|>1-\eta$.

There is also a bound in the opposite direction for the absolute size of the low block. If $d=\dist_{\T}(H,V)>0$, the Dirichlet kernel of $S_{[-L,L]}^\T$ is bounded in modulus by $1/\sin(\pi d)$ on $V-H$. For $g$ normalized as in the lemma,
\[
(1-\eta)p<\int_V|g_0|^2
\le \frac{|V||H|}{\sin^2(\pi d)}.
\]
For our separated intervals below, $|V|=\delta/32$ and $d\ge1/4$. Thus $p=O_\eta(\delta)$ at best, irrespective of the choice of moments. Much smaller values of $p$ are possible. This explains why one cannot simply enlarge the one-step transfer to order one, and why repeating a small spectral selection is the useful alternative. Once $p>0$ is fixed, the required amplification depth is finite by \eqref{eq:amplify}.

\section{Spatial recursion: cells at the first exit}\label{sec:geometry}
Fix $0<\delta<1$. Take
\[
I=[0,1/4],\qquad H\Subset(0,1/4),\qquad C_0=[1/2,5/8],
\]
where $H$ is a nonempty open interval. Choose an open interval $V\Subset C_0^\circ$ of length $(\delta/4)|C_0|$. Apply Lemma~\ref{lem:moment} with these $H,V$ and an accuracy $\eta$ to be fixed below. Extend $g$ periodically to $\R$.

Take $K_1,\ldots,K_N$ to be powers of two, with $K_1\ge16$ and $K_j$ a multiple of $2LK_{j-1}$ for $j>1$. In particular $4K_k\mid K_j$ for $k<j$. We choose the scales successively and may increase each as necessary. For $j\ge1$ define
\begin{equation}\label{eq:fB}
F_j=\{x:\{K_kx\}\in I,\ 1\le k\le j\},
\end{equation}
where braces denote the fractional part. Set $F_0=\R$, $E_0=\varnothing$, and recursively let
\begin{equation}\label{eq:EBj}
A_j=F_{j-1}\cap\{x:\{K_jx\}\in V\},\qquad E_j=E_{j-1}\cup A_j.
\end{equation}
The closed periodic sets $F_j$ are nested, and any product of the first $j$ rescaled factors $g(K_kx)$, with or without modulation, vanishes off $F_j$. A point of $A_j$ satisfies every earlier support condition but fails the $j$-th; no condition is imposed on its later phases. This is exactly the receiving set for the new low-frequency block, since that block will contain the previous product but need not contain any later factors. The sets $A_j$ are pairwise disjoint, and $E_N\cap F_N=\varnothing$.

\begin{lemma}[Thinness independent of the number of steps]\label{lem:geometry}
For all such scales,
\begin{equation}\label{eq:maxdensity}
M\1_{E_N}(y)\le\delta/2\qquad(y\in F_N).
\end{equation}
\end{lemma}
\begin{proof}
At step $j$ consider every cell $C=(m+C_0)/K_j$, $m\in\Z$, meeting $A_j$. Each earlier support condition is constant in the interior of the period cell $(m,m+1)/K_j$: indeed, its boundaries are multiples of $1/(4K_k)$, hence multiples of $1/K_j$ for $k<j$. Since $C$ is contained strictly inside that period cell, it is contained in $F_{j-1}$. Also $C\cap F_j=\varnothing$, since $C_0\cap I=\varnothing$. Thus the selected cells of all steps are pairwise disjoint: the cells at step $j$ lie in $F_{j-1}\setminus F_j$.

The distance from $C_0$ to $I+\Z$ is $1/4$, so
\[
\dist(C,F_N)\ge (4K_j)^{-1}=2|C|.
\]
Moreover $E_N\cap C=A_j\cap C$, and $|E_N\cap C|\le(\delta/4)|C|$. Let $Q$ be an arbitrary bounded interval containing $y\in F_N$. Every selected cell meeting $E_N\cap Q$ has length at most $|Q|/2$. All such cells therefore lie in the concentric interval $2Q$. Disjointness gives the bound for every $Q\ni y$ required by \eqref{eq:maxdensity}:
\[
|E_N\cap Q|\le (\delta/4)\sum_{C:C\cap E_N\cap Q\ne\varnothing}|C|
\le(\delta/4)|2Q|=\tfrac{\delta}{2}|Q|.\qedhere
\]
\end{proof}

\section{Frequency recursion: retain the low block, copy into the high blocks}\label{sec:digits}
The spatial recursion has already fixed the receiving sets. We now arrange the frequency intervals to match them. The guard region in Lemma~\ref{lem:moment} allows us to use dyadic intervals rather than arbitrary endpoints.

Fix $g,L,G,N$. For $D>G$ sufficiently large let
\[
q(t)=\sum_{|n|\le D}a_ne^{2\pi int},\qquad a_n=\widehat g(n),\qquad q_1=q-g_0\ne0.
\]
The only active digits of $q$ are the low digits $|n|\le L/2-1$ and the high digits $G<|n|\le D$. Write $\widetilde g(t)=e^{2\pi iLt}g(t)$ and $\widetilde q(t)=e^{2\pi iLt}q(t)$; the modulation moves the low block into the middle of a dyadic interval.

Work on $X=\R$ or $\T$. Choose a bounded open interval $U$ (an arc on $\T$) and $f_0\in C_c^\infty(U)$ of norm one in $L^2(X)$. Let $P_0=S_{[-R_0,R_0]}^Xf_0$ with $R_0$ large. Define
\[
P_j=P_0\prod_{k=1}^j\widetilde q(K_kx),\qquad u_j=f_0\prod_{k=1}^j\widetilde g(K_kx)\quad(0\le j\le N).
\]
The actual input is $u_j$, supported in $U\cap F_j$; $P_j$ is used only to calculate the projections. Set $B_q=\max(1,\norm q_\infty,\norm g_\infty)$ and
\begin{equation}\label{eq:approx}
e_0=\norm{P_0-f_0}_2 B_q^N+N\norm{q-g}_\infty B_q^{N-1}.
\end{equation}
Telescoping the product difference gives $\norm{P_j-u_j}_{L^2(X)}\le e_0$ for every $j\le N$, uniformly in the scales. With $N$ fixed, choosing $D$ and then $R_0$ large makes $e_0$ arbitrarily small, by uniform Fourier convergence for $g$ and $L^2$ convergence for $f_0$.

At step $j$ choose $K_j$ as in Section~\ref{sec:geometry}, with the additional requirement
\begin{equation}\label{eq:lowerror}
R_{j-1}=R_0+(D+L)\sum_{k<j}K_k<K_j/(16\Lambda).
\end{equation}
The preceding product has its spectrum in $[-R_{j-1},R_{j-1}]$. Multiplication by $\widetilde q(K_jx)$ puts a copy of this band at each active center $(L+n)K_j$. The new central interval collects the low bands; each high band contains a translated copy of the preceding selection. Starting with $\Omega^{(0)}=\varnothing$, set
\begin{equation}\label{eq:intervals}
U_j=[0,2LK_j),\qquad
\Omega^{(j)}=\{U_j\}\ \cup\!\!\bigcup_{G<|n|\le D}\big((L+n)K_j+\Omega^{(j-1)}\big).
\end{equation}
For example, $\Omega^{(1)}=\{U_1\}$, and $\Omega^{(2)}$ consists of $U_2$ and copies of $U_1$ in the high bands of the second factor. This is the entire frequency recursion.

\begin{lemma}[Exact spectral recursion]\label{lem:digits}
Each $\Omega^{(j)}\subset\mathcal D_{\rm dyad}$ has pairwise disjoint closed dilations $\Lambda\omega$ and is contained in $[-R_j,R_j]$. On $X=\R$ or $\T$,
\begin{align}
S_{U_j}^X P_j(x)&=e^{2\pi iLK_jx}g_0(K_jx)P_{j-1}(x),\label{eq:central}\\
S_{(L+n)K_j+\omega}^X P_j(x)&=a_n e^{2\pi i(L+n)K_jx}S_\omega^X P_{j-1}(x)
\quad(G<|n|\le D,\ \omega\in\Omega^{(j-1)}).\label{eq:copy}
\end{align}
Moreover $\norm{P_j}_{L^2(X)}^2=\norm q_{L^2(\T)}^{2j}\norm{P_0}_{L^2(X)}^2$.
\end{lemma}
\begin{proof}
Use induction, writing $R=R_{j-1}$ and $K=K_j$. The translated frequency bands $(L+n)K+[-R,R]$ are disjoint because $2R<K$. Each active low band is within distance $(L/2-1)K+R<LK/2$ of $LK$, and therefore lies inside $U_j$. All high bands are outside $U_j$. Expanding the product in Fourier modes proves \eqref{eq:central}--\eqref{eq:copy}. The same disjoint bands give $\norm{P_j}_2^2=\norm q_2^2\norm{P_{j-1}}_2^2$ by Plancherel.

For the interval geometry, $2LK$ is a power of two, so $U_j$ is dyadic. Every old interval has length $2LK_i$, $i<j$, dividing $K$; translation by $(L+n)K$ preserves its dyadic alignment. The old family lies in $[-R,R]$, and its dilated intervals lie in $[-\Lambda R,\Lambda R]$. Distinct high copies of this enlarged region are disjoint, since $2\Lambda R<K$. The closed dilation of $U_j$ has center $LK$ and radius $\Lambda LK$. A high center is at distance $|n|K\ge(G+1)K$ from $LK$. Since
\[
(G+1)K>\Lambda LK+\Lambda R,
\]
the dilated central interval and every dilated high copy are disjoint, with positive separation. Separation within each copy follows by induction. All undilated intervals lie in $[-R_j,R_j]$ by the definition of $R_j$ (in particular $2LK<(D+L)K$). This completes the proof.
\end{proof}

Define the pointwise squared energy $\mathcal Q_j=\sum_{\omega\in\Omega^{(j)}}|S_\omega^X P_j|^2$, with $\mathcal Q_0=0$, and put $\theta=\norm{q_1}_{L^2(\T)}^2$. Summing \eqref{eq:central}--\eqref{eq:copy} gives the central identity of the construction:
\begin{equation}\label{eq:pointwise-recursion}
\boxed{\mathcal Q_j=|g_0(K_j\cdot)|^2|P_{j-1}|^2+\theta\mathcal Q_{j-1}.}
\end{equation}
The new low block supplies fresh energy, and the high blocks preserve a fraction of the previous energy. These are separate square-function coordinates. This identity retains the cancellations within each original projection while never counting a Fourier coefficient twice.

The smooth version requires no further selection. Fix an even $\psi\in C_c^\infty((-1,1))$ with $0\le\psi\le1$ and $\psi=1$ on $[-1/2,1/2]$. For an interval $\omega$ with center $h$ and radius $r$ set
\begin{equation}\label{eq:smoothsymbol}
m_\omega(\xi)=\psi((\xi-h)/r),\qquad A_\omega^X=m_\omega(D_x).
\end{equation}
At a central step, the active low bands are in the plateau, by the first inequality in the proof of Lemma~\ref{lem:digits}, and the high bands are outside the support. In a high copy the symbol is the translated old symbol. Thus the same projection identities hold for the smooth multipliers by induction, and in particular
\begin{equation}\label{eq:smoothidentity}
A_\omega^X P_N=S_\omega^X P_N\qquad(\omega\in\Omega^{(N)}).
\end{equation}
For arbitrary inputs, both the sharp and these smooth vector-valued multipliers are contractions on unweighted $L^2$, by disjointness and the bound $|m_\omega|\le1$.

\section{Amplification: one scalar recurrence}\label{sec:energy}
Only elementary averaging is needed. For $b$ a bounded one-periodic function and $h\in L^1(X)$, with $X=\R$ or $\T$,
\begin{equation}\label{eq:average}
\int_X h(x)b(Kx)\dd x\longrightarrow \left(\int_Xh\right)\left(\int_0^1 b(t)\dd t\right)
\qquad(K\to\infty,\ K\in\Z).
\end{equation}
For an interval indicator $h$, divide the interval into complete periods and at most two remainders. Approximation by finite step functions in $L^1(X)$ proves the claim; it also holds along the required integer multiples.

We now estimate only the receiving energy. The energy retained in the input region is already controlled by the approximation \eqref{eq:approx} and the exact norm identity in Lemma~\ref{lem:digits}. Set
\[
\sigma=\int_V|g_0|^2,
\qquad b_*=\norm{q}_{L^2(\T)}^{2N}\norm{P_0}_{L^2(X)}^2-e_0^2.
\]
Since $\norm{q}_{L^2(\T)},\norm{P_0}_{L^2(X)}\le1$ and $u_j$ vanishes outside $U\cap F_j$, for $0\le j\le N$ we have
\begin{align*}
\int_{U\cap F_j}|P_j|^2
&=\norm{P_j}_{L^2(X)}^2-\int_{X\setminus(U\cap F_j)}|P_j-u_j|^2\\
&\ge \norm{q}_{L^2(\T)}^{2N}\norm{P_0}_{L^2(X)}^2-e_0^2=b_*.
\end{align*}
In the case $X=\T$, $F_j$ and $E_j$ are understood modulo one. For fixed $N$, $b_*\to1$ as the two approximations are improved, independently of all the chosen scales.

Define
\[
\mathcal H_j=\int_{U\cap E_j}\mathcal Q_j,
\qquad \mathcal H_0=0.
\]
Given $\tau>0$, choose $K_j$ sufficiently large after the preceding scales are fixed. By \eqref{eq:average}, applied to $|P_{j-1}|^2\1_{U\cap F_{j-1}}$ and the periodic function $|g_0|^2\1_V$, the new central block contributes
\begin{equation}\label{eq:newenergy}
\int_{U\cap A_j}|g_0(K_jx)P_{j-1}(x)|^2\dd x
\ge \sigma\int_{U\cap F_{j-1}}|P_{j-1}|^2-\tau
\ge \sigma b_*-\tau.
\end{equation}
Integrate the first term of \eqref{eq:pointwise-recursion} on $U\cap A_j$ and the second on $U\cap E_{j-1}$. Both sets are contained in $U\cap E_j$; both terms are nonnegative. By \eqref{eq:newenergy},
\begin{equation}\label{eq:renewal}
\mathcal H_j\ge \theta\mathcal H_{j-1}+\sigma b_*-\tau.
\end{equation}
By induction,
\begin{equation}\label{eq:limitenergy}
\mathcal H_N\ge (\sigma b_*-\tau)\sum_{j=1}^N\theta^{N-j}.
\end{equation}
As $D,R_0$ are taken large and $\tau$ small, the right side approaches
\begin{equation}\label{eq:energybound}
\mathcal L_N=\sigma\sum_{j=1}^N(1-p)^{N-j}>(1-\eta)[1-(1-p)^N],
\end{equation}
because $\theta\to1-p$. This is the amplification budget \eqref{eq:amplify}. The exact same calculation holds for the smooth symbols on $P_N$ by \eqref{eq:smoothidentity}.

\begin{lemma}[Concentration inside an arbitrary prescribed window]\label{lem:periodic}
For $X=\R$ or $\T$, any bounded open interval or arc $U\subset X$, and $0<\delta,\kappa<1$, there are measurable $E,F\subset U$, $u\in C_c^\infty(U)\setminus\{0\}$ with $\supp u\subset F$, and a finite family $\Omega\subset\mathcal D_{\rm dyad}$ with pairwise disjoint closed dilations $\Lambda\omega$, such that
\begin{equation}\label{eq:periodic}
M_X\1_E\le\delta\ \text{on }F,\qquad
\sum_{\omega\in\Omega}\int_E |S_\omega^X u|^2>\kappa\norm u_{L^2(X)}^2.
\end{equation}
The second inequality also holds simultaneously for the smooth multipliers $A_\omega^X$. Here $M_\T$ is the maximal operator over arcs.
\end{lemma}
\begin{proof}
Choose $\rho\in(\sqrt\kappa,1)$ and $\eta$ with $1-\eta>\rho^2$. First obtain $g,p$ from Lemma~\ref{lem:moment}, then choose $N$ so that $\mathcal L_N>\rho^2$. With $N$ fixed, take $D$ and then $R_0$ sufficiently large that
\[
\sigma b_*\sum_{j=1}^N\theta^{N-j}>\rho^2,
\qquad e_0<e:=(\rho-\sqrt\kappa)/3.
\]
Choose $\tau>0$ so small that the right side of \eqref{eq:limitenergy} still exceeds $\rho^2$. Finally choose $K_1,\ldots,K_N$ successively to satisfy \eqref{eq:newenergy}, \eqref{eq:lowerror} and the divisibility conditions. Then $\mathcal H_N>\rho^2$.

Take $E=E_N\cap U$, $F=F_N\cap U$, $u=u_N$ and $\Omega=\Omega^{(N)}$. The function $u$ is smooth with compact support in $U$ and its support is contained in $F_N$. Lemma~\ref{lem:geometry} proves thinness, by restriction on the line or by periodically extending $E$ on the circle (the extension is contained in the original periodic $E_N$). The spectral recursion gives $\norm{P_N}_{L^2(X)}\le1$. Therefore $\norm u_2\le1+e$. By the contraction property and \eqref{eq:approx}, the square root of the receiving energy exceeds $\rho-e>\sqrt\kappa(1+e)$. This proves \eqref{eq:periodic} and nonvanishing. The smooth inequality follows in exactly the same way from \eqref{eq:smoothidentity}.
\end{proof}

\begin{proof}[Proof of Theorem~\ref{thm:thin}]
Apply Lemma~\ref{lem:periodic} with $X=\R$ and any bounded open interval $U$. Plancherel gives the upper bound one for the supremum; the lemma gives the matching lower supremum. We will also use the simultaneous smooth conclusion, which we record as
\begin{equation}\label{eq:smooththin}
\sum_{\omega\in\Omega}\int_E |A_\omega^\R u|^2>\kappa\norm u_2^2.
\end{equation}
\end{proof}

The arbitrary window $U$ will be useful when constructing a single periodic weight.

\section{A single counterexample weight and interval family}\label{sec:weight}
We recall the Coifman--Rochberg maximal-weight fact \cite{CR80}: if $0<\alpha<1$ and $E$ is a measurable set of positive measure (bounded or periodic in the applications below), then
\begin{equation}\label{eq:CR}
v_E=(M\1_E)^\alpha\in A_1,\qquad [v_E]_{A_1}\le C_\alpha,
\end{equation}
where $C_\alpha$ is independent of $E$. Here is a short proof of the version we need. For a bounded interval $K$ and $x\in K$, split $\1_E$ into its restrictions to $3K$ and $(3K)^c$, where $3K$ is the concentric triple. For $y\in K$, any interval containing $y$ and reaching $(3K)^c$ has length at least $\abs K$; its enlargement to include $x$ has length at most twice the original length. Thus $M(\1_{E\setminus3K})(y)\le2M\1_E(x)$. The weak $(1,1)$ inequality for $M$ applied to $u_K=\1_{E\cap3K}$, followed by integration of the distribution function, gives
\[
\frac1{\abs K}\int_K(Mu_K)^\alpha\lesssim_\alpha
\left(\frac{\norm{u_K}_1}{\abs K}\right)^\alpha
\le 3^\alpha (M\1_E(x))^\alpha.
\]
In the first step, integrate $\alpha\lambda^{\alpha-1}\min\{\abs K,C\norm{u_K}_1/\lambda\}$ over $\lambda>0$ and divide by $\abs K$. Subadditivity of the $\alpha$-power now bounds the average of $v_E$ on $K$ by a constant times $v_E(x)$. Taking the supremum over $K\ni x$ gives \eqref{eq:CR}.

For every $0<\alpha<1$, Theorem~\ref{thm:thin} with $\kappa=1/2$ gives finite examples for which $v_E=1$ almost everywhere on $E$ and $v_E\le\delta^\alpha$ on $F$. Thus
\begin{equation}\label{eq:varying}
\frac{\norm{T_\Omega f}_{L^2(v_E)}^2}{\norm f_{L^2(v_E)}^2}>\frac{1}{2\delta^\alpha}.
\end{equation}
For $0<t\le1$, Jensen's inequality gives $M(v^t)\le(Mv)^t$ and hence $[v^t]_{A_1}\le[v]_{A_1}^t$. In particular, weights of the form $(M\1_E)^\alpha$ can have uniformly prescribed characteristic $1+\varepsilon$ by taking $\alpha>0$ small enough. Rather than repeating the construction for each $\varepsilon$, we now fix one weight and one family.

\begin{proof}[Proof of Theorem~\ref{thm:weighted}]
\emph{First fix the frequency family.} Take $\delta_n=2^{-n}$ and $\kappa_n=1-2^{-n}$. Choose finite examples $(E_n,F_n,f_n,\Omega_n)$ by Theorem~\ref{thm:thin}, with these parameters and the fixed dilation factor $\Lambda$. Use the simultaneous smooth conclusion \eqref{eq:smooththin} as well.

For each $n$, choose a positive dyadic scale $s_n=2^{-m_n}$ sufficiently small that the affine images $c_n+s_n(\Lambda\omega)$, $\omega\in\Omega_n$, are contained in the interior of $J_n=[2^{-n-1},2^{-n})$, where $c_n$ is the midpoint of $J_n$. We may also arrange that $c_n$ is an integer multiple of every scaled interval length $s_n|\omega|$: the lengths and $c_n$ are dyadic, and the family is finite. The transformed intervals $c_n+s_n\omega$ then remain in the standard grid. With
\[
f'_n(x)=s_n^{1/2}e^{2\pi i c_nx}f_n(s_nx),\qquad
E'_n=s_n^{-1}E_n,\qquad F'_n=s_n^{-1}F_n,
\]
maximal density and the unweighted receiving fraction are unchanged. Indeed, dilation and modulation give the exact projection identity
\[
S_{c_n+s_n\omega}f'_n(x)=s_n^{1/2}e^{2\pi i c_nx}(S_\omega f_n)(s_nx),
\]
and the same formula holds for the corresponding rescaled smooth symbols. We relabel the transformed data as $(E_n,F_n,f_n,\Omega_n)$. Now
\begin{equation}\label{eq:fixedfamily}
\Omega_*=\bigcup_{n\ge1}\Omega_n\subset\mathcal D_{\rm dyad},\qquad
\bigcup_{\omega\in\Omega_*}\omega\subset(0,1).
\end{equation}
The closed $\Lambda$-dilations are pairwise disjoint, since different subfamilies lie, even after dilation, strictly inside distinct $J_n$.

\emph{Then fix the spatial weight.} Set $\alpha=1/2$, $b_n=2^{-n}$, and $v_n=(M\1_{E_n})^\alpha$. These satisfy $0<v_n\le1$, \eqref{eq:CR}, and $v_n(x)\to0$ as $|x|\to\infty$: for $x$ far from bounded $E_n$, every interval containing $x$ and meeting $E_n$ has length at least $\dist(x,E_n)$.

For a translation $t_n>0$ set
\[
w_n(x)=\tfrac12\big(v_n(x-t_n)+v_n(-x-t_n)\big).
\]
Each $w_n$ is even, bounded by one, decays at infinity, and has the same uniform $A_1$ bound. Choose the translations inductively so that
\begin{align}
\sup_{x\in F_n+t_n}\sum_{k<n}b_kw_k(x)&\le \tfrac12 b_n\delta_n^\alpha,\label{eq:past}\\
\sup_{x\in F_k+t_k}b_nw_n(x)&\le2^{-(n-k+1)}b_k\delta_k^\alpha\quad(k<n),\label{eq:future}\\
\sup_{x\in F_n+t_n}v_n(-x-t_n)&\le\delta_n^\alpha.\label{eq:selfreflection}
\end{align}
There are finitely many conditions at each step, all attainable by taking $t_n$ large: all input sets are bounded, the earlier profiles decay at infinity, and the last condition samples $v_n$ on $-F_n-2t_n$.

Set
\begin{equation}\label{eq:W}
W=\sum_{n\ge1}b_nw_n.
\end{equation}
The weight is even and $0<W\le1$. Sublinearity and monotone convergence give $MW\le C_{1/2}W$. On $F_n+t_n$, its own term is at most $b_n\delta_n^\alpha$ by thinness and \eqref{eq:selfreflection}, and the earlier and later terms contribute at most $b_n\delta_n^\alpha$ in total by \eqref{eq:past}--\eqref{eq:future}. On $E_n+t_n$, the original translated component is $b_n/2$ almost everywhere. Hence
\begin{equation}\label{eq:budget}
W\le2b_n\delta_n^\alpha\ \text{on }F_n+t_n,
\qquad W\ge b_n/2\ \text{a.e. on }E_n+t_n.
\end{equation}

\emph{Use every positive power.} Let $u_n(x)=f_n(x-t_n)$. The unweighted receiving fraction for $u_n$ and $\Omega_n$ exceeds $\kappa_n$, by translation invariance. Since $\Omega_n\subset\Omega_*$, \eqref{eq:budget} gives, for every fixed $t>0$,
\begin{equation}\label{eq:fixedweight}
\frac{\norm{T_{\Omega_*}u_n}_{L^2(W^t)}^2}{\norm{u_n}_{L^2(W^t)}^2}
>\frac{\kappa_n(b_n/2)^t}{(2b_n\delta_n^\alpha)^t}
=\kappa_n4^{-t}\delta_n^{-t/2}\longrightarrow\infty.
\end{equation}
All the test functions are smooth and compactly supported and all the numerators are finite, by $W^t\le1$ and \eqref{eq:plancherel}. For $0<t\le1$, Jensen gives $[W^t]_{A_1}\le[W]_{A_1}^t$ as above. Letting $t\downarrow0$ proves the near-one characteristic assertion with the same $\Omega_*$. No assertion that $W^t\in A_1$ for every $t>1$ is needed. The same argument works with the selected smooth symbols. This completes the proof.
\end{proof}

\begin{remark}[One BMO direction]\label{rem:bmo}
Put $B=-\log W$. The $A_1$ condition ensures a positive essential infimum of $W$ on every bounded interval $I$. If $m_I=\operatorname*{ess\,inf}_I W$, then Jensen gives $|I|^{-1}\int_I\log(W/m_I)\le\log[W]_{A_1}$. Since $\log(W/m_I)\ge0$, its mean oscillation is at most twice its mean. Thus $B\in\mathrm{BMO}$, and the same fixed family has norm one at $t=0$ and infinite norm on $L^2(e^{-tB})$ for every $t>0$. The perturbation $tB$ tends to zero in BMO; this is not a claim of uniform pointwise convergence of $W^t$ to one.
\end{remark}

\section{Sharp contrast and the quadratic boundary}\label{sec:thresholds}
Two short consequences describe the size of the failure and the sequence exponent at which it occurs.

\begin{corollary}[Exact contrast saturation]\label{cor:contrast}
Fix the family $\Omega_*$ from Theorem~\ref{thm:weighted}. For every $0<a<b$ and $\varepsilon>0$,
\begin{equation}\label{eq:contrast}
\sup_{\substack{w\ \mathrm{even},\ a\le w\le b\\ [w]_{A_1}\le1+\varepsilon}}
\norm{T_{\Omega_*}}_{L^2(w)\to L^2(w)}=\sqrt{b/a}.
\end{equation}
The supremum is unchanged if the weights are required to be smooth. The same value is obtained by taking the supremum over weights and finite subfamilies of $\Omega_*$.
\end{corollary}
\begin{proof}
The upper bound is immediate from $a\le w\le b$ and Plancherel. For the lower bound, use the examples of the preceding proof, with $\kappa_n\to1$. Choose $0<t\le1$ so small that $[W^t]_{A_1}\le1+\varepsilon$ and set
\[
\nu_n=a+(b-a)\min\{1,(2W/b_n)^t\}.
\]
The $A_1$ characteristic is preserved under scaling and truncation from above: if $Mv\le C v$, then $M\min(1,v)\le\min(1,Mv)\le C\min(1,v)$. Positive sums preserve the common bound. Thus $\nu_n$ is even, $a\le\nu_n\le b$, and $[\nu_n]_{A_1}\le1+\varepsilon$. By \eqref{eq:budget}, $\nu_n=b$ on $E_n+t_n$ almost everywhere and $\nu_n\le a+(b-a)(4\delta_n^{1/2})^t$ on the input set. Consequently
\[
\frac{\norm{T_{\Omega_n}u_n}_{L^2(\nu_n)}^2}{\norm{u_n}_{L^2(\nu_n)}^2}
>\frac{b\kappa_n}{a+(b-a)(4\delta_n^{1/2})^t}\longrightarrow b/a.
\]
This proves both assertions about the supremum. Finally convolve each $\nu_n$ with an even, nonnegative smooth compactly supported probability kernel tending to the identity. The pointwise bounds and $A_1$ bound are preserved, since maximal averaging commutes with translation and $M(\nu_n*\phi)\le(M\nu_n)*\phi$. Dominated convergence preserves the two norm integrals for each fixed $u_n$ and finite family, by the bounded weights and Plancherel. This proves the smooth assertion.
\end{proof}

\begin{corollary}[The quadratic endpoint of Kr\'ol's estimates]\label{cor:q}
For $1\le s<\infty$ set
\[
T_{\Omega,s}f=\left(\sum_{\omega\in\Omega}|S_\omega f|^s\right)^{1/s},
\]
with the supremum when $s=\infty$. For $2<s<\infty$, Kr\'ol's Theorem~B(i) \cite{Krol14} gives
\begin{equation}\label{eq:krol-positive}
\norm{T_{\Omega,s}f}_{L^2(w)}\lesssim_{s,w}\norm f_{L^2(w)},\qquad w\in A_{2-2/s},
\end{equation}
uniformly over disjoint interval families. The case $s=\infty$, $w\in A_2$, follows from weighted Carleson estimates. At the limiting exponent $s=2$, the same assertion fails for $A_1$ by Theorem~\ref{thm:weighted}. Hence the uniform strong estimate for every $A_1$ weight holds precisely for $2<s\le\infty$. The fixed family and weights $W^t$, $0<t\le1$, already give failure for $1\le s\le2$; for $s<2$ there is failure even without weights.
\end{corollary}
\begin{proof}
In Kr\'ol's notation the sequence exponent is $q'=q/(q-1)$. Set $q=s/(s-1)$ and take spatial exponent $p=2$ in his Theorem~B(i). Then $p/q=2-2/s$, proving the known positive estimate \eqref{eq:krol-positive}. For the supremum over intervals, write an interval projection as the difference of two modulated Hilbert transforms and apply the weighted Carleson bound; see Hunt--Young for Fourier series \cite{HY74} and the real-line estimates in \cite{DL14}.

Theorem~\ref{thm:weighted} and $T_{\Omega,s}\ge T_{\Omega,2}$ give failure for $s\le2$. For the unweighted claim take $N$ consecutive unit intervals and $\widehat f=\1_{[0,N]}$. Each component has the same modulus as the inverse transform of the unit-interval indicator. Plancherel gives norm ratio $N^{1/s-1/2}$, which diverges for $s<2$. Finally, for the fixed family and any fixed $W^t$, $0<t\le1$, the constants in the bounds for $s>2$ must diverge as $s\downarrow2$: apply the bounds to each finite subfamily, take $s\downarrow2$, and then exhaust the family by monotone convergence.
\end{proof}

For completeness, the weaker positive assertion with $w\in A_1$ also follows by a direct interpolation argument. The weighted Carleson bound gives $L^{p_-}(w)\to L^{p_-}(w;\ell^\infty)$ for the linear vector of projections when $1<p_-<2$, and Rubio de Francia's theorem \cite{RF85} gives $L^{p_+}(w)\to L^{p_+}(w;\ell^2)$ for $p_+>2$. For $s>2$, put $\vartheta=2/s$ and choose $p_+>2$ sufficiently close to two that $1/2=(1-\vartheta)/p_-+\vartheta/p_+$ has a solution $1<p_-<2$. Complex interpolation on the common measure space $(\R,w\,\dd x)$ gives the $L^2(w;\ell^s)$ estimate. This is an alternative proof only in the $A_1$ case; the stronger weight range in \eqref{eq:krol-positive} is due to Kr\'ol. The index $2-2/s$ decreases to $1$ as $s\downarrow2$. We do not assert a necessary weight-class characterization for each individual $s>2$.

\section{Smooth cutoffs and maximal-function bumps}\label{sec:smooth}
For a fixed positive integer $s$, define the intrinsic smooth square function, in the Fourier normalization of this paper, by
\[
G_\Omega^{(s)} f(x)=\left(\sum_{\omega\in\Omega}
\sup_{\phi\in\Phi_\omega^{(s)}}\abs{\mathcal F^{-1}(\phi\widehat f)(x)}^2\right)^{1/2},
\quad
\Phi_\omega^{(s)}=\{\phi\in C_c^\infty(\omega^\circ):
\max_{0\le k\le s}\abs\omega^k\norm{\phi^{(k)}}_\infty\le1\}.
\]
Up to the Fourier convention this is the intrinsic operator of \cite[p.~447]{DFPR25}; taking compact support strictly inside each interval only restricts the supremum.

\begin{corollary}[Smooth interval multipliers]\label{cor:smooth}
For the same weight and family in Theorem~\ref{thm:weighted}, for every fixed $s$ and every $t>0$, $G_{\Omega_*}^{(s)}$ is unbounded on $L^2(W^t)$. In fact, the square function of the particular smooth multipliers in \eqref{eq:smoothsymbol}, rescaled and shifted along with their intervals, is unbounded after multiplication by a fixed positive constant meeting the normalization of $\Phi_\omega^{(s)}$.
\end{corollary}
\begin{proof}
The proof of Theorem~\ref{thm:weighted} uses the simultaneous examples \eqref{eq:thin} and \eqref{eq:smooththin}. The affine frequency maps and spatial translations preserve the smooth selection identities and receiving fractions; the same budgets \eqref{eq:budget} therefore give \eqref{eq:fixedweight} for the smooth square function and every $t>0$. For an interval of radius $r$, the symbol $\psi((\xi-h)/r)$ satisfies
\[
\abs\omega^k\norm{m_\omega^{(k)}}_\infty=2^k\norm{\psi^{(k)}}_\infty.
\]
Hence multiplication by $a_s=(\max_{0\le k\le s}2^k\norm{\psi^{(k)}}_\infty)^{-1}>0$ puts all the symbols into $\Phi_\omega^{(s)}$, uniformly in the interval family. The intrinsic square function dominates the resulting square function pointwise. Since $a_s$ is independent of $n$, its ratios still diverge. The same choice also meets the derivative adaptation bounds on p.~446 of \cite{DFPR25}.
\end{proof}

\begin{corollary}[Failure of finite maximal bumps]\label{cor:bumps}
For each $0<r<\infty$ and each integer $k\ge1$, neither of the estimates
\[
\int T_\Omega f^2w\lesssim_r\int\abs f^2(M(w^r))^{1/r},
\qquad
\int T_\Omega f^2w\lesssim_k\int\abs f^2 M^kw
\]
holds uniformly over disjoint interval families and positive weights, even if the $A_1$ characteristics of the weights are uniformly bounded. Here $M^k$ denotes $k$ successive applications of $M$.
\end{corollary}
\begin{proof}
Choose $0<\alpha<\min(1,1/r)$ and use the weights $v_E=(M\1_E)^\alpha$ from \eqref{eq:varying}. By \eqref{eq:CR}, $M(v_E^r)\le C_{\alpha r}v_E^r$, and $M^k v_E\le C_\alpha^k v_E$. Both proposed input weights are therefore bounded above by a constant times $v_E$. The ratio in \eqref{eq:varying} tends to infinity as $\delta\downarrow0$. The same proof works for the smooth multipliers by \eqref{eq:smooththin}.
\end{proof}

The preceding failure contrasts with the classical theory of Calder\'on--Zygmund singular integrals. In the power-bump setting initiated by C\'ordoba and Fefferman \cite{CF76}, such operators admit a bound with input weight $(M(w^r))^{1/r}$ for each $r>1$. P\'erez \cite{Perez94} gives the stronger iterated-maximal estimate at exponent two, with input weight $M^3w$. Indeed the power-bump bound also follows from that estimate: $(M(w^r))^{1/r}$ is an $A_1$ majorant of $w$, by the Coifman--Rochberg maximal-weight theorem, and hence dominates $M^3w$ up to a constant depending on $r$. For $0<r\le1$ the power bump is no larger than $Mw$, by Jensen's inequality. The failure above for every finite iteration and every positive power therefore distinguishes arbitrary-frequency square functions from individual classical singular integrals.

\section{A periodic even-weight consequence}\label{sec:periodic-even}
The periodic counterpart with even weights and real Fourier coefficients merits separate consideration in view of \cite{Berkson22}. We give a fixed-weight version, without imposing a uniform dependence on the characteristic in the proposed inequality.

\begin{corollary}\label{cor:periodic-even}
For each $\varepsilon>0$ and fixed $\Lambda\ge1$ there exist an even weight $W_{\T}\in A_1(\T)$, with $[W_{\T}]_{A_1(\T)}\le1+\varepsilon$, and a countable dyadic family with integer endpoints and pairwise disjoint closed dilations $\Lambda\omega$ such that the periodic Rubio de Francia square function is unbounded on $L^2(\T,W_{\T})$, even on smooth functions all of whose Fourier coefficients are real. The maximal operator defining $A_1(\T)$ is taken over arcs.
\end{corollary}
\begin{proof}
Choose points $s_n\in(1/8,1/4)$ and open arcs $U_n$ containing them such that the closures of these arcs and their reflections are mutually disjoint. Set $\delta_n=2^{-n}$. For any small subarc $U\subset U_n$ centered at $s_n$, Lemma~\ref{lem:periodic} gives a finite periodic example $(E_n,F_n,f_n,\Omega_n)$ supported in $U$ with density at most $\delta_n/2$ (use the stronger bound in its proof) and receiving fraction greater than $1/2$. Define
\[
v_n=(M_\R\1_{E_n+\Z})^\alpha,\qquad w_n(x)=\tfrac12(v_n(x)+v_n(-x)).
\]
In fact $M_\R\1_{E_n+\Z}\le\delta_n/2$ on $F_n$: the periodically extended $E_n$ is contained in the first-exit sets of Lemma~\ref{lem:geometry}. The maximal-weight proof \eqref{eq:CR} applies to these periodic indicators and gives a common $A_1(\R)$ bound, hence also the same bound on $\T$. Choose $0<\alpha<1$ sufficiently small, using \eqref{eq:CR} and the Jensen power inequality of Section~\ref{sec:weight}, to give characteristic at most $1+\varepsilon$.

As $|U|\to0$, $v_n\to0$ uniformly away from $s_n$ on the circle, regardless of the example within $U$. Indeed, an interval of length at most $1/2$ meeting $E_n+\Z$ and a point a fixed distance away must have length bounded below by that distance, and its receiving measure is at most $|U|$. For longer intervals the average is at most $5|U|$, by periodicity. Since $0<v_n\le1$, dominated convergence gives $\int_\T w_n\to0$. The reflected contribution is also uniformly small on $F_n$, so by choosing $U$ sufficiently short we have $w_n\le\delta_n^\alpha$ on $F_n$ and $w_n\ge1/2$ almost everywhere on $E_n$.

Choose positive coefficients $c_n$ inductively so large that $\sum_{k<n}c_k\le c_n\delta_n^\alpha/2$. After choosing $c_n$, choose the arc and its example as above so that additionally
\[
c_n\int_\T w_n\le2^{-n},\qquad
\sup_{F_k}c_n w_n\le2^{-(n-k+1)}c_k\delta_k^\alpha\quad(k<n).
\]
The disjoint arcs and uniform decay off $\pm s_n$ make these finitely many requirements possible. Then $W_{\T}=\sum_n c_n w_n$ is even, integrable, positive almost everywhere, and has the common $A_1(\T)$ bound by monotone convergence. The same budgets as \eqref{eq:budget} hold with $c_n,F_n,E_n$ in place of $b_n,F_n+t_n,E_n+t_n$. Therefore the squared norm ratios of these finite periodic examples are greater than $\tfrac18\delta_n^{-\alpha}$. Translate each finite family by a sufficiently large integer multiple of its largest dyadic length, so that its closed dilations are disjoint from the earlier dilated families. This preserves dyadic alignment and translates the projections by modulation of the test functions. The resulting family is fixed. No bounded frequency window is asserted on the circle.

Finally, define the conjugate-linear involution $\mathcal J u(x)=\overline{u(-x)}$. It is an isometry on $L^2(\T,W_{\T})$ because the weight is even, and commutes with every real-symbol Fourier multiplier, including the projections. Both $u_1=(u+\mathcal J u)/2$ and $u_2=(u-\mathcal J u)/(2i)$ are fixed by $\mathcal J$, so their Fourier coefficients are real. The parallelogram identity gives $\norm{u_1}^2+\norm{u_2}^2=\norm u^2$ in the weighted input space, and the same identity for the sum of the squared projected output norms, applied to any finite subfamily. At least one of the two smooth functions therefore has a ratio no smaller than that of $u$. Applying this to the diverging sequence completes the proof. The same conjugate-reflection decomposition applies to the real-line examples with the even weight $W$ in Theorem~\ref{thm:weighted}; hence the real-line unboundedness also holds on smooth compactly supported inputs whose Fourier transforms are real-valued.
\end{proof}

\section{Relation to earlier results}\label{sec:relation}
\subsection{Weak versus strong control and quadratic localization}
There is no conflict with weighted weak-type bounds in \cite{Krol14,GRS21,DFPR25}. To see what such bounds still imply for the constructed configurations, put $v_E=(M\1_E)^\alpha$. For fixed $\alpha$, the weak endpoint gives
\[
\sup_{\lambda>0}\lambda^2\abs{\{x\in E:T_\Omega f(x)>\lambda\}}
\lesssim_\alpha\delta^\alpha\norm f_2^2.
\]
On the other hand, Theorem~\ref{thm:thin} makes the integrated receiving energy an order-one fraction of $\norm f_2^2$. A bound on the weak distribution functional does not bound its integral with respect to $\dd\lambda/\lambda$. The construction permits an arbitrarily large range of amplitudes and supplies no uniform higher-integrability bound as the number of factors grows.

A related difficulty arises in sparse localization by quadratic averages. For example, put $I_j=[0,2^{-j}]$, $0\le j\le N$, and $f=\1_{I_N}$. Then
\[
\int_{\R}\sum_{j=0}^N\1_{I_j}(x)\frac{1}{|I_j|}\int_{I_j}|f(y)|^2\dd y\dd x
=(N+1)\norm f_2^2.
\]
The example concerns the positive averages, not the Fourier operator: repeated averaging loses the endpoint energy budget. Ordinary interpolation from below two also fails, since for $\widehat f=\1_{[0,N]}$ and $N$ unit frequency intervals the $L^p$ norm ratio grows as $N^{1/p-1/2}$ when $1<p<2$.

\subsection{Square functions versus scalar variation multipliers}
Remark~5.4, pp.~30--31, of Lacey's monograph \cite{Lacey07} (Remark~5.61 in its arXiv version) states the strong $A_1$ square-function inequality as a valid estimate. Theorem~\ref{thm:weighted} contradicts that formula as printed, when the constant is independent of the disjoint intervals. No proof of the endpoint formula is supplied in the remark.

There are two distinct issues. First, even granting the sharp-function estimate (5.3) for the smooth square function $G^\Omega$ in the preceding Remark~5.2, the usual weighted sharp-function inequality leads only to
\[
\norm{G^\Omega f}_{L^2(w)}^2\lesssim_w \int M(|f|^2)\,w.
\]
To conclude the claimed endpoint from this estimate one would need strong $L^1(w)$ control of $M$. This is false even for $w=1$: the maximal function of $\1_{[0,1]}$ has a nonintegrable tail. The $A_1$ property controls $Mw$, not the strong $L^1(w)$ norm of $M$. This explains why the preceding sharp-function estimate is insufficient; we do not claim that Lacey actually used this inference to obtain Remark~5.4.

Second, the moment lemma reveals a specific error in the proof of the relative off-support estimate (2.19) on pp.~12--14: the step on p.~13 invokes rapid spatial decay to estimate the near part of an already projected input relative to its own norm. The support condition is on the original input, not on its projection. Appendix~\ref{app:literature} locates this step explicitly and gives a density argument disproving the estimate, as well as the analogous (3.6) in Proposition~3.3. These observations concern particular estimates and the printed weighted formula, not the classical unweighted theorem. The original formulation in \cite[Section~6, (6.4)]{RF85} and the explicit later open-problem discussion in \cite[Section~1.3]{DFPR25} describe the same endpoint.

The adjacent variation-multiplier statement must be distinguished from the square-function statement. Xu's Theorem~A and his principal higher-dimensional theorem \cite{Xu96} concern unweighted $L^p$. The weighted atomic statement $R_2(\mathcal D)\subset M_2(\R,w)$ for $w\in A_1$ is recorded in Kr\'ol \cite[Theorem~3]{Krol14}, with attribution to Th\'eor\`eme~1 and Lemme~5 of \cite{CRS88}. The weighted variation result $V_2(\mathcal D)\subset M_2(\R,w)$ is contained in the Berkson--Gillespie theorem as formulated in \cite[Theorem~4]{Krol14} (attributed to \cite[Theorem~1.2]{BG98}) and also in Kr\'ol's Theorem~A(i). Here $\mathcal D$ denotes the classical lacunary dyadic frequency partition, not the full grid $\mathcal D_{\rm dyad}$. Kr\'ol's Remark~9 sketches the unweighted variation argument using the classical Littlewood--Paley and Rubio de Francia inequalities and interpolation; we have not inspected the original four-page note \cite{CRS88} directly. We do not infer a verification of its weighted proof from that sketch. In any event, the scalar bounds just stated are compatible with our counterexample, for the following reason.

For clarity, fix a finite family $\Omega$ and $f\in L^2(w)\cap L^2$. The map $A_f:\ell^2(\Omega)\to L^2(w)$, $A_f z=\sum_\omega z_\omega S_\omega f$, satisfies
\[
\norm{A_f}_{\mathrm{op}}\lesssim_w\norm f_{L^2(w)},\qquad
\norm{A_f}_{\mathrm{HS}}^2=\sum_\omega\norm{S_\omega f}_{L^2(w)}^2.
\]
The first estimate follows from the scalar theorem because the step multiplier $\sum_\omega z_\omega\1_\omega$ has uniformly bounded dyadic $2$-variation norm when $\norm z_{\ell^2}\le1$: on any ordered partition each constant segment can contribute to at most two nonzero increments, and $|u-v|^2\le2(|u|^2+|v|^2)$. Endpoints may be assigned half-open conventions without changing the multiplier. The second quantity, however, is the Hilbert--Schmidt norm, not the operator norm. It can grow with rank even when the first norm is bounded. In particular, the random sign vectors in the introductory identity have $\ell^2$ norm $|\Omega|^{1/2}$, not one. This is why scalar weighted variation bounds and failure of the strong square-function bound are compatible.

\subsection{Congruent subdivisions of lacunary intervals}
Garg, Roncal and Shrivastava \cite[Section~4.2]{GRS21} obtain a positive endpoint case by first partitioning the frequency line into lacunary intervals $J_k$ and then subdividing each $J_k$ into congruent intervals, with lengths allowed to depend on $k$. Their estimate follows by first applying the lacunary square-function bound to the parent pieces and then the congruent-interval bound within each parent. The recursion \eqref{eq:intervals} has a different geometry: a single high block contains a translated multiscale family, rather than a congruent subdivision. Its members can have very different lengths while lying within a tiny relative frequency neighborhood. The two-stage estimate of \cite{GRS21} is therefore not a uniform estimate for this arbitrary-depth recursion. Applying a separate bound at every recursive step would introduce depth-dependent constants. In particular, membership of the intervals in the standard dyadic grid and separation of their prescribed dilations must not be confused with the classical decomposition into lacunary annuli.

\subsection{Earlier affirmative claims}
Theorem~\ref{thm:weighted} uses an even weight, but not a radially decreasing one, and is compatible with the even decreasing-weight results of \cite{Berkson23,DFPR25}. Corollary~\ref{cor:periodic-even} addresses the periodic same-weight assertions in Theorems~4 and 5 of \cite{Berkson22}. Appendix~\ref{app:literature} records the relevant displays from that paper and gives independent tests of its intermediate periodic estimate and its real-line estimates. The radial proof in \cite{DFPR25} uses the sufficient covering condition (43), p.~468, to obtain the wave-packet bound (42): a single centered interval suffices at each level. The Fourier counterexample does not impose that radial structure.

\subsection{The Walsh model}
The distinction from the Walsh model is structural. Lemma~\ref{lem:moment} relies on the fact that all members of a finite set of trigonometric characters remain linearly independent upon restriction to any nonempty open interval. This does not hold for Walsh characters on a dyadic cell. For example, the constant Walsh character and the first Rademacher function coincide on $[0,1/2)$. Their moments cannot be prescribed independently by functions supported there. More generally, projection onto the first $2^m$ Walsh characters is dyadic conditional expectation: it is constant on each dyadic atom $Q$ of length $2^{-m}$. Thus if $V$ occupies at most a fraction $\rho$ of each such atom, then $\int_V\abs{\mathbb E_m g}^2\le\rho\norm{\mathbb E_m g}_2^2$. This rules out the corresponding concentration mechanism at dyadic low-pass scales. The actual proof of the Walsh endpoint in \cite[Theorem D and Section~5]{DFPR25} uses exact localized orthogonality of Walsh packets, giving maximal-weight control by the dyadic maximal operator. The Fourier construction here does not preserve that locality.

\subsection{Scope of the obstruction}\label{sec:scope}
The one-dimensional example has consequences for other frequency geometries. We record two elementary reductions, in the context of higher-dimensional and directional square functions (see, for example, \cite{ADHPR23}). Neither requires a new amplification argument.

For $d\ge2$, put $\widetilde W(x)=W(x_1)$. For axis-parallel cubes, its $A_1$ characteristic equals the one-dimensional characteristic of $W$. With the fixed family $\Omega_*$ of Theorem~\ref{thm:weighted}, take the dyadic cubes $Q_\omega=\omega\times[0,|\omega|)^{d-1}$. They retain disjoint closed $\Lambda$-dilations. Given a finite subfamily, choose a Schwartz function $h$ of unit $L^2$ norm on $\R^{d-1}$ whose Fourier support is in $(0,\min_\omega|\omega|)^{d-1}$. Then $S_{Q_\omega}(f\otimes h)=(S_\omega f)\otimes h$ for every interval in the subfamily. The weighted norm ratios are identical, so the square function of this single countable cube family is unbounded on $L^2(\widetilde W^t)$ for all $t>0$. For ball-based $A_1$ one obtains an equivalent characteristic bound; small powers still tend to characteristic one.

In the plane, the intervals in $(0,1)$ also define pairwise disjoint conical sectors $C_\omega=\{(\xi_1,\xi_2):\xi_2>0,\ \xi_1/\xi_2\in\omega\}$. For $\widehat h\in C_c^\infty((-1,1))$ of unit $L^2$ norm, put $\widehat h_\epsilon(\xi_2)=\epsilon^{-1/2}\widehat h((\xi_2-1)/\epsilon)$ with $0<\epsilon<1/2$. For each fixed $f\in L^2(\R)$ and each finite family,
\[
\sum_\omega\norm{S_{C_\omega}(f\otimes h_\epsilon)-(S_\omega f)\otimes h_\epsilon}_2^2\longrightarrow0
\quad(\epsilon\downarrow0).
\]
Indeed, after the substitution $\xi_2=1+\epsilon s$, the multiplier indicators converge almost everywhere in $(\xi_1,s)$, and dominated convergence and Plancherel apply. Since $\widetilde W^t\le1$, the same convergence holds in the weighted norm. The input norm factors exactly; using the divergent finite-subfamily tests of Theorem~\ref{thm:weighted} and choosing $\epsilon$ separately for each test gives unboundedness for the fixed sector family. In both reductions, smooth compactly supported tests follow by unweighted $L^2$ approximation and the disjoint-projection contraction. These are particular higher-dimensional endpoint obstructions, not objections to known directional estimates with different weight or geometry hypotheses.

The failure of finite maximal bumps in Corollary~\ref{cor:bumps} also excludes any uniform domination whose right-hand side has the corresponding strong weighted bound. This does not contradict the quadratic-average sparse localization in \cite{DFPR25}, which need not give such an endpoint bound.

The distinction between absolute spectral energy and relative spatial concentration is central to the present construction. Its combination with recursive frequency selection provides a mechanism for producing endpoint obstructions beyond the original weighted conjecture, as the extensions to dyadic cubes and conical sectors illustrate. The next subsection establishes a local transfer to elliptic spectral projections. Multilinear square functions and endpoint sparse-domination problems remain further directions.

\subsection{Local transfer to elliptic spectral projections}\label{sec:elliptic}
The obstruction also passes from Fourier intervals to spectral intervals of elliptic operators. This is a local statement: it is unaffected by the boundary condition or by the behaviour of the coefficients far from the test functions. We give a fixed-weight version, rather than only a failure of uniformity over weights.

Let $D\subset\R^d$ be open and let $A$ be the nonnegative self-adjoint Dirichlet or Neumann realization in $L^2(D)$ of
\[
A=-\nabla\cdot(a(x)\nabla)+V(x),
\]
where $a$ is bounded, real symmetric and uniformly positive definite, and $V\ge0$ is locally bounded. The realizations are defined by their closed quadratic forms. Assume that $a$ is $C^1$ in an interior ball $B_0\Subset D$. For a family $\mathcal I$ of disjoint positive spectral intervals write
\[
\mathcal S_{A,\mathcal I}f=\left(\sum_{J\in\mathcal I}|\1_J(A)f|^2\right)^{1/2}.
\]
Here the spectral projections are those of the unweighted self-adjoint operator. We use the uncentered maximal operator $M_d$ over axis-parallel cubes for ambient weights on $\R^d$.

\begin{theorem}[Local spectral obstruction]\label{thm:elliptic}
Under these assumptions there exist an ambient weight $w\ge1$ in $A_1(\R^d)$ and a fixed countable family $\mathcal I_*$ of pairwise disjoint bounded intervals in $(0,\infty)$ such that
\begin{equation}\label{eq:elliptic-unbounded}
\sup_{0\ne f\in C_c^\infty(B_0)}
\frac{\norm{\mathcal S_{A,\mathcal I_*}f}_{L^2(D,w^t)}}{\norm f_{L^2(D,w^t)}}=\infty
\qquad(t>0).
\end{equation}
Every positive power $w^t$ is locally integrable, and $[w^t]_{A_1}\le[w]_{A_1}^t$ for $0<t\le1$. In particular a counterexample weight of characteristic at most $1+\varepsilon$ is available for every $\varepsilon>0$. The intervals can be grouped in finite families whose spectral locations escape to infinity.
\end{theorem}

\begin{proof}
We first transfer the thin-set theorem to a constant-coefficient operator, then freeze the coefficients at an interior point, and finally assemble one weight.

\emph{1. Constant-coefficient spectral intervals.} Let $A_0=-\nabla\cdot(a_0\nabla)$ on $\R^d$, where $a_0$ is a fixed real positive definite matrix. Given $0<\delta,\kappa<1$, we claim there are bounded $E,F\subset\R^d$, $0\ne u\in C_c^\infty(\R^d)$ supported in $F$, and a finite family $\mathcal I$ of disjoint bounded positive spectral intervals such that
\begin{equation}\label{eq:tangent-thin}
M_d\1_E\le\delta\ \text{on }F,\qquad
\int_E\sum_{J\in\mathcal I}|\1_J(A_0)u|^2>\kappa\norm u_2^2.
\end{equation}
Start with a one-dimensional configuration of Theorem~\ref{thm:thin} with a receiving fraction strictly greater than $\kappa$. Modulate its input $f$ by $e^{2\pi icx_1}$ with $c$ large enough that all translated intervals $c+\omega$ lie in $(0,\infty)$. For each $\omega$ let $J_\omega=4\pi^2(a_0)_{11}(c+\omega)^2$, the image of the interval under the increasing positive quadratic map. On the axis $\xi'=0$ the multiplier $\1_{J_\omega}(4\pi^2\xi^Ta_0\xi)$ selects both $c+\omega$ and $-(c+\omega)$. The former gives the desired modulated Fourier output; the total squared $L^2$ error from the latter is bounded by the integral of $|\widehat f|^2$ over the disjoint intervals $-2c-\omega$, which tends to zero as $c\to\infty$. Thus for sufficiently large $c$ the one-dimensional receiving fraction for these spectral projections still exceeds $\kappa$.

If $d>1$, use $u=f_c\otimes h_\eta$, where $f_c=e^{2\pi ic\cdot}f$ and $h_\eta$ is a normalized Schwartz function whose Fourier support shrinks to $0$ in the transverse variables. For instance take $h_\eta(x')=\eta^{(d-1)/2}h(\eta x')$ with $\widehat h\in C_c^\infty$ and $\norm h_2=1$. Substituting $\xi'=\eta s$ in the Plancherel formula shows by dominated convergence that the vector of spectral outputs approaches the tensor product of the one-dimensional outputs with $h_\eta$ in the direct-sum $L^2$ norm. Indeed the symbol $4\pi^2\xi^Ta_0\xi$ tends to $4\pi^2(a_0)_{11}\xi_1^2$, and the finitely many boundary frequencies have measure zero. Fix $\eta$ small enough to preserve a strict margin above $\kappa$, restrict the receiving set to $E_1\times[-R,R]^{d-1}$, and truncate $h_\eta$ smoothly to compact support. Taking $R$ and the truncation radius sufficiently large preserves the margin, by monotone convergence, $L^2$ approximation and the contraction property of the disjoint spectral projections. Here $E_1,F_1$ are the original one-dimensional sets. The resulting input is supported in $F_1$ times a bounded transverse box. For every cube $Q$ containing an input point $y$, the average of $\1_E$ over $Q$ is at most the average of $\1_{E_1}$ over its first-coordinate interval. Hence $M_d\1_E(y)\le M\1_{E_1}(y_1)\le\delta$. This proves \eqref{eq:tangent-thin}; for $d=1$ the transverse step is unnecessary.

\emph{2. Freezing coefficients.} Fix $x_0\in B_0$, put $D_h=(D-x_0)/h$, and let $U_hu(x)=h^{-d/2}u((x-x_0)/h)$, a unitary map from $L^2(D_h)$ to $L^2(D)$. On $L^2(\R^d)=L^2(D_h)\oplus L^2(\R^d\setminus D_h)$ extend $h^2U_h^{-1}AU_h$ by the zero operator on the second summand; call the resulting nonnegative self-adjoint operator $A_h$. For $\phi\in C_c^\infty(\R^d)$, its support is in the smooth interior region after scaling for all sufficiently small $h$. It is in the domain of $A_h$, and
\[
A_h\phi=-\nabla\cdot(a(x_0+h\cdot)\nabla\phi)+h^2V(x_0+h\cdot)\phi\longrightarrow A_0\phi
\quad\text{in }L^2,
\]
where $A_0=-\nabla\cdot(a(x_0)\nabla)$. This uses only local coefficient regularity. The set $(A_0+1)C_c^\infty(\R^d)$ is dense in $L^2$: by Fourier transformation the domain of $A_0$ is $H^2$, and $C_c^\infty$ is dense in its graph norm. On this dense set the resolvent identity gives
\[
(A_h+1)^{-1}(A_0+1)\phi-\phi=(A_h+1)^{-1}(A_0-A_h)\phi\longrightarrow0.
\]
Uniform resolvent bounds extend convergence to all of $L^2$. Thus $A_h\to A_0$ in the strong resolvent sense. The standard spectral convergence theorem \cite[Theorem~VIII.24]{RS80} yields $\1_J(A_h)u\to\1_J(A_0)u$ strongly for every bounded interval $J$ whose endpoints are not eigenvalues of $A_0$. The constant-coefficient operator has no spectral atoms, since its quadratic symbol has Lebesgue-null level sets. Half-open intervals are allowed by the same spectral convergence at both endpoints.

Apply \eqref{eq:tangent-thin} to this $A_0$ with a receiving fraction strictly greater than a prescribed $\kappa<1$. The finite collection of output norms on $E$ converges as $h\downarrow0$. Scaling back, for all sufficiently small $h$ there are $E_h=x_0+hE$, $F_h=x_0+hF$, $f_h=U_hu$, supported in $F_h$, and the disjoint spectral family $h^{-2}\mathcal I$ with
\begin{equation}\label{eq:local-spectral-thin}
M_d\1_{E_h}\le\delta\ \text{on }F_h,\qquad
\int_{E_h}\sum_{J\in h^{-2}\mathcal I}|\1_J(A)f_h|^2>\kappa\norm{f_h}_2^2.
\end{equation}
The density condition is invariant under the same spatial dilation and translation. The sets and input support fit into any prescribed interior neighbourhood of $x_0$ for sufficiently small $h$. Since the original spectral intervals are bounded away from zero, we may at the same time put all the rescaled spectral intervals above any prescribed finite energy. This argument does not depend on the choice of boundary condition.

\emph{3. One ambient weight and one spectral family.} Choose a sequence of interior balls $B_n$ with mutually disjoint closures contained in $B_0$, and centers $x_n$. Fix $0<\alpha<1$. Inductively put $b_n=n^2(1+\sum_{k<n}b_k)$ and choose $\delta_n$ so small that $b_n\delta_n^\alpha\le1$. Use the construction in step 2 at $x_n$ with these parameters, $\kappa=1/2$, and $h_n$ small enough that both $E_n,F_n\subset B_n$. Denote the resulting input and spectral family by $f_n,\mathcal I_n$. Further shrink $h_n$ so that $\mathcal I_n$ lies above all earlier spectral intervals. Set $v_n=(M_d\1_{E_n})^\alpha$. The multidimensional Coifman--Rochberg theorem \cite{CR80} gives $[v_n]_{A_1}\le C_{d,\alpha}$, independently of the sets. Also $0<v_n\le1$, $v_n=1$ almost everywhere on $E_n$, and $v_n\le\delta_n^\alpha$ on $F_n$.

For each fixed stage $n$, the unscaled receiving set is bounded. As $h_n\downarrow0$, $v_n\to0$ uniformly on sets a positive distance from $x_n$ and pointwise away from $x_n$. This follows by bounding the measure of $E_n$ by $O(h_n^d)$ and noting that any cube meeting $E_n$ and a distant point must have side length bounded below. Thus we can impose the additional finitely many requirements
\[
\sup_{F_k}b_nv_n\le2^{-n}\quad(k<n),\qquad
\sup_{\R^d\setminus B_0}b_nv_n\le2^{-n},
\qquad
\norm{b_nv_n}_{L^n([-n,n]^d)}\le2^{-n}.
\]
The last follows from dominated convergence at the fixed stage. All requirements are compatible with \eqref{eq:local-spectral-thin}, which holds for all sufficiently small scales. Define $w=1+\sum_n b_nv_n$ and $\mathcal I_*=\bigcup_n\mathcal I_n$. The local norm conditions give $w\in L^m_{\rm loc}$ for each finite integer $m$ by H\"older and summation. Moreover $w\le2$ outside $B_0$. The common $A_1$ bound passes to this positive sum by monotone convergence. Write $s_{n-1}=1+\sum_{k<n}b_k$. On $F_n$ we have $w\le s_{n-1}+1+\sum_{k>n}2^{-k}\le3s_{n-1}$, while $w\ge b_n$ almost everywhere on $E_n$. Therefore, for every $t>0$,
\[
\frac{\norm{\mathcal S_{A,\mathcal I_*}f_n}_{L^2(D,w^t)}^2}{\norm{f_n}_{L^2(D,w^t)}^2}
\ge\frac12\left(\frac{b_n}{3s_{n-1}}\right)^t
=\frac12(n^2/3)^t\longrightarrow\infty.
\]
The finite spectral outputs belong to the form domain and hence to $H^1(D)$; their restrictions to $B_0$ belong to some $L^{2+\epsilon_d}(B_0)$ by Sobolev embedding. Since $w$ has every finite local integrability exponent and is bounded outside $B_0$, each finite-family weighted output norm is finite, by H\"older. Jensen's inequality gives the asserted characteristic bound for small powers. This proves the theorem.
\end{proof}

\noindent\emph{Consequences.} The theorem covers the Laplacian, uniformly elliptic constant-coefficient operators, and locally $C^1$ variable-coefficient divergence-form operators with a nonnegative locally bounded potential, on the whole space or with the stated boundary conditions. It also gives a stationary PDE benchmark: if $u=(1+A)^{-1}f$, then $(1+A)\1_J(A)u=\1_J(A)f$. Consequently the sum of squared weighted norms of these frequency-resolved elliptic residuals cannot be bounded by $C_{A,w}\norm f_{L^2(D,w)}^2$ for arbitrary spectral intervals. This is an obstruction to an unrestricted spectral refinement of regularity, not to ordinary elliptic regularity or to standard dyadic spectral multiplier estimates. We have not imposed a dyadic grid or prescribed dilation separation in the spectral variable of the general operator. The corresponding problem on Riemannian manifolds requires a coordinate and weight transfer, which is not pursued here.

\appendix
\section{Quantitative constraints and amplification cost}\label{app:cost}
The counterexample needs no quantitative lower estimate for $p$. For completeness, we record what is known about the minimum complexity and the cost of the particular guarded construction used here. These estimates are not needed in the main proof.

\subsection{General necessary complexity bounds}
The following estimates are independent of the counterexample. Throughout, $E,F$ are bounded, $|E|>0$, $0\ne f\in L^2(\R)$ is supported in $F$, $M\1_E\le\delta$ on $F$, and the frequency intervals are disjoint. We write $\mathcal E=\int_E T_\Omega f^2$ for the receiving energy. Finite families suffice for the proofs; monotone convergence handles countable ones.

\begin{lemma}[Necessary complexity]\label{lem:complexity}
The following bounds hold with absolute constants, except for the indicated dependence on $\alpha$.
\begin{enumerate}
\item If the lengths belong to $N$ classes comparable within a factor of two, then, for $0<\delta<e^{-2}$,
\begin{equation}\label{eq:universal-depth}
\mathcal E\lesssim\min\{1,N\delta\log^3(1/\delta)\}\norm f_2^2.
\end{equation}
\item If $F$ is closed and the receiving set meets only $J$ dyadic distance layers $\{2^j\le\dist(x,F)<2^{j+1}\}$ up to null sets, then $\mathcal E\lesssim\min\{1,J\delta\}\norm f_2^2$.
\item If $E,F$ lie in $K$ ordered consecutive blocks, each containing points of at most one of the two sets, then $\mathcal E\lesssim\min\{1,\delta(1+\log K)^2\}\norm f_2^2$.
\item For $f$ bounded and any $0<\alpha<1$,
\begin{equation}\label{eq:amplitudes}
\mathcal E\lesssim_\alpha\delta^\alpha\norm f_2^2
\left[1+\log_+\frac{|E|\norm f_\infty^2}{\delta^\alpha\norm f_2^2}\right].
\end{equation}
\end{enumerate}
\end{lemma}
\begin{proof}
\emph{Length classes.} Scale a single class to lengths in $[1,2]$. The centers $c_k$ are at least one apart. Take $\phi\in C_c^\infty(\R)$ equal to one on $[-1,1]$, and put $K_\phi=\check\phi$. The one-dimensional Sobolev inequality on the disjoint unit cells about $c_k$, followed by Plancherel in $y$, applied to the Fourier transform of $f(x-y)K_\phi(y)$ gives
\[
\sum_k |(\phi(\cdot-c_k)\widehat f)^\vee(x)|^2
\lesssim\int|f(x-y)|^2(1+y^2)|K_\phi(y)|^2\dd y.
\]
The rapidly decreasing nonnegative kernel on the right convolved with a weight $w$ is bounded by $C Mw$ (split into dyadic annuli). Integration therefore bounds the smooth squared norm by $C[w]_{A_1}\norm f_{L^2(w)}^2$. On each smooth piece recover $S_\omega$ by a difference of two modulated half-line projections. Petermichl's sharp Hilbert-transform bound \cite{Pet07,Hyt12} gives their $L^2(w)$ norm at most $C[w]_{A_2}$, and $[w]_{A_2}\le[w]_{A_1}$. Thus the sharp squared norm for one class is at most $C[w]_{A_1}^3\norm f_{L^2(w)}^2$. The proof of \eqref{eq:CR} gives $[(M\1_E)^\alpha]_{A_1}\lesssim(1-\alpha)^{-1}$ for $1/2\le\alpha<1$, by integrating the weak $(1,1)$ bound. With this weight, summing over the $N$ classes gives $CN(1-\alpha)^{-3}\delta^\alpha\norm f_2^2$. Choose $\alpha=1-1/\log(1/\delta)$ and also use Plancherel.

\emph{Distance layers.} We use a simple sampling inequality: if $A$ occupies at most $\rho d$ in each cell of a grid of length $d$, then
\begin{equation}\label{eq:appendix-sampling}
\int_A|h|^2\lesssim\rho(\norm h_2^2+d^2\norm{h'}_2^2),\qquad h\in H^1(\R).
\end{equation}
Indeed, on each cell $I$, $\sup_I|h|^2\lesssim d^{-1}\int_I|h|^2+d\int_I|h'|^2$; multiply by $|A\cap I|$ and sum. After modulation, Bernstein's inequality gives the bound $C\rho\norm h_2^2$ for frequency support of length at most $d^{-1}$.

In particular the entire near-field contribution, where $D(x)|\omega|\le B$, $D(x)=\dist(x,F)$, satisfies
\begin{equation}\label{eq:appendix-near}
\sum_\omega\int_{E\cap\{D(x)|\omega|\le B\}}|S_\omega f(x)|^2\dd x\lesssim B\delta\norm f_2^2\qquad(B\ge1).
\end{equation}
For each $\omega$, use cells of length $|\omega|^{-1}$; a cell meeting this set can be enlarged to an interval of length $O(B|\omega|^{-1})$ containing a point of $F$. Its receiving measure is therefore $O(B\delta|\omega|^{-1})$. Apply \eqref{eq:appendix-sampling} and sum by Plancherel. The closure of $F$ may be used here, since the uncentered maximal function is lower semicontinuous.

Fix a distance layer $E_d=E\cap\{d\le D(x)<2d\}$. Each cell of length $d$ meeting $E_d$ has receiving measure $O(\delta d)$. Short intervals $|\omega|\le d^{-1}$ are paid by sampling. For the long intervals take a smooth even cutoff $\chi$ supported in $[-1,1]$, equal to one near zero, and set
\[
K_d(z)=\frac{1-\chi(z/d)}{2\pi iz},\qquad h_{a,d}=K_d*(e^{-2\pi ia\cdot}f).
\]
On $E_d$, for $\omega=[a,b]$ the exact kernel formula is $S_\omega f=e^{2\pi ibx}h_{b,d}-e^{2\pi iax}h_{a,d}$. The multiplier $\widehat K_d(\xi)=m(d\xi)$ is bounded and rapidly decreasing at infinity: boundedness follows from the sine integral for the truncated odd Hilbert kernel and decay from integration by parts (every positive-order derivative of $K_1$ is integrable). The left endpoints of the long intervals are $d^{-1}$-separated, as are the right endpoints separately. Consequently Plancherel gives, for either endpoint family,
\[
\sum_a(\norm{h_{a,d}}_2^2+d^2\norm{h'_{a,d}}_2^2)
=\int|\widehat f(\xi)|^2\sum_a(1+4\pi^2d^2|\xi-a|^2)|m(d(\xi-a))|^2\dd\xi
\lesssim\norm f_2^2.
\]
Apply \eqref{eq:appendix-sampling} to these envelopes on $E_d$. This proves $\int_{E_d}T_\Omega f^2\lesssim\delta\norm f_2^2$. Summing over $J$ layers proves part (2); almost every point of $E$ is outside $F$ by thinness and differentiation.

\emph{Ordered blocks.} First suppose $F\subset(-\infty,s_0]$, $E\subset[s_0,\infty)$. Translate the rightmost point of the closure of $F$ to zero, so the input is supported in $(-\infty,0]$ and $M\1_E(0)\le\delta$. Then $D(x)=x$ for $x>0$. Split $E$ into $E_j=E\cap[2^j,2^{j+1})$, $d_j=2^j$. The contribution $|\omega|d_j\le A$ for a fixed sufficiently large $A$ is paid by \eqref{eq:appendix-near} with $B=2A$. For each $j$ cover $E_j$ by the grid cells $I$ of length $d_j/8$ that meet it. Their union $U_j$ is contained in $[7d_j/8,17d_j/8]$, so the $U_j$ have bounded overlap. The receiving measure in each cell is $O(\delta d_j)$ by the interval reaching zero.

For $x>0$ put $h_a(x)=(2\pi i)^{-1}\int_{-\infty}^0 e^{-2\pi iay}f(y)/(x-y)\dd y$. For $|b-a|d_j>A$, the two modulations of the constant means of $h_a,h_b$ on $I$ are uniformly separated in $L^2(I)$ (integrate the cross exponential). Poincar\'e's inequality for the deviations from the means yields
\[
\int_I(|h_a|^2+|h_b|^2)\lesssim\int_I |S_{[a,b]}f|^2+d_j^2\int_I(|h_a'|^2+|h_b'|^2).
\]
Local Sobolev sampling and the endpoint kernel identity now bound the receiving energy in that cell by $C\delta$ times the right side. Sum the original projection term using bounded overlap and Plancherel. It remains to bound, before multiplication by $C\delta$,
\begin{equation}\label{eq:appendix-res}
\sum_\omega\sum_{a\in\partial\omega}\int_{x\ge cA/|\omega|}x^2|h_a'(x)|^2\dd x.
\end{equation}
For the one-sided input define the Hardy-space function $G(a+it)=\int_{-\infty}^0 f(y)e^{-2\pi iay}e^{2\pi ty}\dd y$ for $t>0$. Its $H^2$ norm is bounded by $\norm f_2$ by Plancherel. The Laplace representation $(x-y)^{-2}=4\pi^2\int_0^\infty t e^{-2\pi t(x-y)}\dd t$ and Cauchy--Schwarz give
\[
x^2|h_a'(x)|^2\lesssim\int_0^\infty t e^{-2\pi tx}|G(a+it)|^2\dd t.
\]
Integrating in $x$ in \eqref{eq:appendix-res} produces the measure $\mu=\sum_\omega\sum_{a\in\partial\omega}\delta_a\otimes e^{-c'A t/|\omega|}\dd t$. For any interval $J_0$ of length $\ell$,
\[
\mu(J_0\times(0,\ell))\lesssim\sum_\omega\sum_{a\in\partial\omega\cap J_0}\min(\ell,|\omega|)\lesssim\ell.
\]
Indeed, the intervals inside $J_0$ have total length at most $\ell$, and only a bounded number cross its boundary. This is the Carleson box condition. The classical $H^2$ embedding gives $\int |G|^2\dd\mu\lesssim\norm f_2^2$: the nontangential maximal function is bounded on $L^2$ by the maximal function of the boundary modulus; its open level sets decompose into intervals, and the boxes over those intervals cover the corresponding superlevel set in the half-plane. Integrating the distribution bound proves the embedding. Thus \eqref{eq:appendix-res} is paid, giving $\int_E T_\Omega f^2\lesssim\delta\norm f_2^2$ for one-sided transfer. Reflection gives the reverse orientation.

For $K$ ordered blocks build a balanced binary tree of consecutive groups. At each node retain the two transfers between opposite children. Each has operator norm $O(\sqrt\delta)$ by the one-sided bound. At a fixed generation the input and output sets of distinct nodes are disjoint, so the entire generation has the same bound. Each source--receiver pair belongs to exactly one such operator, at the lowest common ancestor. The triangle inequality over at most $O(1+\log K)$ generations proves part (3).

\emph{Amplitude.} Apply the known weighted weak endpoint \cite{Krol14,GRS21,DFPR25} and the weighted strong $L^4$ theorem \cite{RF85} to $v_E=(M\1_E)^\alpha$. Writing $A_0=\delta^\alpha\norm f_2^2$, $B_0=\norm f_\infty$, the distribution on $E$ of $T_\Omega f$ is bounded by the minimum of $|E|$, $C_\alpha A_0/\lambda^2$ and $C_\alpha A_0B_0^2/\lambda^4$. Integrate against $2\lambda\dd\lambda$, splitting at $(A_0/|E|)^{1/2}$ and $B_0$. If the order of these thresholds reverses, the measure and fourth-moment estimates alone give $O_\alpha(A_0)$. This proves \eqref{eq:amplitudes}.
\end{proof}

For a receiving fraction bounded below by $\kappa>0$, these inequalities require respectively $N\gtrsim_\kappa 1/(\delta\log^3(1/\delta))$, $J\gtrsim_\kappa\delta^{-1}$, $1+\log K\gtrsim_\kappa\delta^{-1/2}$, and $1+\log_+(|E|\norm f_\infty^2/(\delta^\alpha\norm f_2^2))\gtrsim_{\kappa,\alpha}\delta^{-\alpha}$. These are necessary estimates, not sharp universal thresholds. The recursion permits all these complexity parameters to grow without bound. The specific guarded mechanism has a sharper constraint, as follows.

\subsection{Cost of the guarded one-factor construction}
Fix $H\Subset(0,1/4)$, $0<\eta<1$ and $\Lambda\ge1$. Let $V$ be centered at $v=9/16$ with length $\delta/32$. For the normalized factors of Lemma~\ref{lem:moment}, the zero moments extend from $L/2-1$ to $G=\lceil(\Lambda+2)L\rceil$. If $p_*(\delta,\eta,\Lambda)$ is the supremum of $p$ over these factors and admissible powers of two $L$, then there are positive constants $c,C$ depending on $H,\eta,\Lambda$ such that, for sufficiently small $\delta$,
\begin{equation}\label{eq:expdepth}
 e^{-C/\delta}\le p_*(\delta,\eta,\Lambda)\le e^{-c/\delta}.
\end{equation}
This determines the exponential order, not an exact constant in the exponent.

For the upper bound write $g_0=\sqrt p\,b$, where $\norm b_2=1$ and $b$ has degree $l\le L/2-1$. Then $(2l+1)|V|\ge1-\eta$. Put $m=G-l$ and $r(t)=\cos^2(\pi(t-v))$. We have $r\le1/2$ on $H$ and $r\ge r_*>1/2$ on $V$, with $r_*$ independent of $\delta$. The polynomial $r^m b$ has degree at most $G$; the guard moments imply $S_{[-G,G]}^\T g=g_0$, so
\[
\sqrt p(1-\eta)r_*^m\le|\langle g_0,r^m b\rangle|
=|\langle g,r^m b\rangle|\le2^{-m}.
\]
Thus $p\le(1-\eta)^{-2}(2r_*)^{-2m}$, with $m\gtrsim_\eta\delta^{-1}$.

For the lower bound choose the normalized Dirichlet polynomial in Lemma~\ref{lem:moment}, with $l\asymp_\eta\delta^{-1}$, $L$ the smallest power of two at least $2l+2$, and $G=\lceil(\Lambda+2)L\rceil\asymp_{\eta,\Lambda}\delta^{-1}$. Use the following elementary norming estimate for trigonometric polynomials $P$ of degree at most $m_0\ge1$ on any fixed subinterval $J\Subset H$:
\begin{equation}\label{eq:norming}
\norm P_{L^2(\T)}\le e^{C_J m_0}\norm P_{L^2(J)}.
\end{equation}
To see this, put $n=2m_0$ and write $e^{2\pi i m_0t}P(t)=Q(e^{2\pi it})$ with $\deg Q\le n$. In $J$ choose $n+1$ alternating subintervals of length $|J|/(2(n+1))$, separated by gaps of that length. In each choose a point $t_j$ where $|P(t_j)|\le (2(n+1)/|J|)^{1/2}\norm P_{L^2(J)}$. Then for $j\ne k$, $|e^{2\pi it_j}-e^{2\pi it_k}|\ge 2|J||j-k|/(n+1)$, since $|J|<1/2$. Lagrange interpolation on the unit circle gives
\[
\sup_{|z|=1}|Q(z)|\le \max_j |P(t_j)|
\left(\frac{n+1}{|J|}\right)^n\frac{2^n}{n!}.
\]
The factorial bound $n!\ge(n/e)^n$ proves \eqref{eq:norming}.

Now fix a nonnegative $\phi\in C_c^\infty(H)$ bounded below by a positive constant on $J$. The moment Gram matrix $\mathsf G_{mn}=\int\phi(t)e^{2\pi i(n-m)t}\dd t$, $|m|,|n|\le G$, satisfies $\mathsf G\ge e^{-C_H G}\operatorname{Id}$ by \eqref{eq:norming} with $m_0=G$. Extend the coefficients of $b$ by zero to degree $G$ and solve $\mathsf Gc=\widehat b$. With $u=\phi\sum c_ne^{2\pi int}$,
\[
\norm u_2^2\le\norm\phi_\infty c^*\mathsf Gc
\le\norm\phi_\infty e^{C_H G}\norm b_2^2.
\]
Then $g=u/\norm u_2$ has the required moments, and $p=\norm u_2^{-2}\ge e^{-C/\delta}$ after adjusting constants. This proves \eqref{eq:expdepth}.

For a fixed desired fraction $0<\kappa<1-\eta$, the ideal receiving budget in \eqref{eq:amplify} is obtained at depth of order $p^{-1}$, with constants depending on the chosen margin. The smallest depth of this guarded one-factor amplification, at fixed $H,\eta,\kappa,\Lambda$, is therefore $\exp(\Theta(1/\delta))$ as $\delta\downarrow0$. The lower bound refers to amplification by the specified geometric-series budget, not to the total possible receiving energy from every interval at every receiver. It is not a universal depth lower bound for arbitrary disjoint families, nor does it give a sharp exponent constant. The real-line and periodic constructions permit the matching upper exponential order through the above explicit choice of moments.

\section{Details of the comparison with earlier claims}\label{app:literature}
The arguments below do not use Theorems~\ref{thm:weighted} or \ref{thm:thin}. They isolate specific intermediate estimates or explicit constants in the published texts; they are not needed to prove the main counterexample.

\subsection{Off-support spectral data and relative localization}
For a bounded measurable set $B\subset\R^d$ and a nonempty open set $U\subset\R^d$,
\begin{equation}\label{eq:remote-density}
\{\widehat f|_B:f\in C_c^\infty(U)\}\quad\text{is dense in }L^2(B).
\end{equation}
Indeed, extend any annihilator $b\in L^2(B)$ by zero. Its inverse Fourier transform is entire (the integral is over a bounded set and $b\in L^1$). Orthogonality to the displayed range implies that the inverse transform vanishes on $U$. Real-analytic uniqueness and Plancherel imply $b=0$. This proves \eqref{eq:remote-density}. Equivalently, for the sharp projection $P_B=\1_B(D)$ and any measurable receiving set $A$, provided $|B|>0$,
\[
\sup_{\substack{f\in C_c^\infty(U)\\P_B f\ne0}}
\frac{\norm{\1_A P_B f}_2}{\norm{P_B f}_2}
=\norm{\1_A P_B}_{L^2\to L^2}.
\]
Indeed the normalized nonzero projections from $U$ are dense in the unit sphere of the bandlimited subspace. Relative concentration of the projected component therefore does not depend on the location of the permitted open input set. This is not a bound on the full input cost.

In particular, if a nonzero wave packet $\psi_s$ and a bounded smoothing symbol $m$ satisfy $\supp\widehat{\psi_s}\cup\supp m\subset B$, there are $f_n\in C_c^\infty(U)$ with $\widehat f_n|_B\to\widehat{\psi_s}|_B$. Then $\langle f_n,\psi_s\rangle\to\norm{\psi_s}_2^2$ and $m(D)f_n\to m(D)\psi_s$ in $L^2$.

For example, in the notation of Lacey \cite[Proposition~2.16, p.~12]{Lacey07}, the asserted estimate (2.19) is
\[
\sum_{s\in\mathcal T(\{\omega\}),\ I_s\subset I}|\langle f,\psi_s\rangle|^2
\lesssim (t\rho)^{-5}\norm{\psi^{3\omega}*f}_2^2,
\qquad \supp f\subset(tI)^c.
\]
Fix $I,\omega$ and a nonzero packet with $I_s\subset I$, with $I$ large enough for the stated condition on $\rho$. For each large $t$, choose $U$ to be an open interval outside $tI$ and $B$ containing both $\supp\widehat{\psi_s}$ and the support of the smoothing multiplier $\widehat{\psi^{3\omega}}$. Applying \eqref{eq:remote-density} and passing to the limit at each fixed $t$ forces
\[
\norm{\psi_s}_2^4\lesssim (t\rho)^{-5}\norm{\psi^{3\omega}*\psi_s}_2^2.
\]
Letting $t\to\infty$ is a contradiction.

The problematic inference in the printed proof can be located more precisely. On p.~13, put $P f=\psi^{3\omega}*f$. While decomposing $P^2f$, the text defines $F_0=\1_{(t/2)I}Pf$ and, using rapid decay, asserts $\norm{F_0}_2\lesssim(t\rho)^{-10}\norm{Pf}_2$. Neither this assertion nor the version with $P^2f$ in place of $Pf$ follows from the support condition on $f$. In fact, choose a fixed bandlimited $h$ with $P^k h\ne0$ ($k=1$ or $2$). By \eqref{eq:remote-density}, for each $t$ one may approximate its band data by an input supported outside $tI$. The asserted bound would then give $\norm{\1_{(t/2)I}P^k h}_2\lesssim(t\rho)^{-10}\norm{Ph}_2$. As $t\to\infty$ the left side tends to $\norm{P^k h}_2>0$. This shows that the relative decay inference fails, independently of the number of convolutions in the definition of $F_0$.

Similarly, Proposition~3.3, p.~20, states (3.6) with $\1_{[-1,1]^d}\le\widehat\varphi\le\1_{[-2,2]^d}$:
\[
\sum_{R_s\subset U_0,\ \omega_s=\omega}|\langle f,\varphi_s\rangle|^2
\lesssim a^{15d}\norm{\varphi^{3\omega}*f}_2^2,
\qquad \supp f\subset\{M_{\rm str}\1_{U_0}>a\}^c.
\]
Here $M_{\rm str}$ is the strong maximal function and $0<a<1$. Fix a bounded rectangle $U_0$ containing a tile rectangle $R_s$; for each $a$ choose a nonempty open box far enough from $U_0$ in every coordinate that $M_{\rm str}\1_{U_0}\le a$ there. Such a box exists because for a rectangle $U_0$ the strong maximal average of its indicator is the product of the one-dimensional maximal averages. Use this box as $U$ in \eqref{eq:remote-density}, with a fixed bounded frequency set containing the packet and smoothing supports. The asserted inequality gives $\norm{\varphi_s}_2^4\lesssim a^{15d}\norm{\varphi^{3\omega}*\varphi_s}_2^2$, again impossible as $a\downarrow0$. On p.~21, the corresponding proof step asserts smallness of $F_0=\1_{\{M_{\rm str}\1_{U_0}>2a\}}(\varphi^{3\omega}*f)$ relative to $\norm{\varphi^{3\omega}*f}_2$ by rapid decay and the support condition on $f$. The same density argument shows why this inference is invalid. Both obstructions concern the relative estimate in terms of the already projected input, not the usual full-input spatial decay estimates (such as (2.23)), or the classical unweighted Rubio de Francia theorem, which has other established proofs.

\subsection{Berkson's 2022 periodic argument}
Theorem~4 of \cite{Berkson22}, pp.~12--15, asserts the periodic square-function bound with constant $2\sqrt C$ for even $A_1(\T)$ weights of characteristic at most $C$ and inputs with real Fourier coefficients. Theorem~5, p.~16, removes the coefficient restriction with constant $4\sqrt C$. Theorem~3 instead retains a length parameter. We reproduce the relevant equalities and estimates below in the angular parametrization $t\in[-\pi,\pi]$, with $e_q(t)=e^{iqt}$ and normalized Fourier coefficients. This makes the objections independent of access to the journal.

Let $F^{[+]}$ be the sum of the positive-coefficient Fourier terms of a trigonometric polynomial $F$ with real Fourier coefficients, and write $S_1(F)=\{q:\widehat F(q)>0\}$. The last equality in the display immediately after (2.17), p.~13, is the assertion
\[
\left\|\left(\sum_k\left|\sum_{q\in S_1(F)\cap V_k}\widehat F^{[+]}(q)e_q\right|^2\right)^{1/2}\right\|_{L^2(\T,w)}
=\left(\int_0^{2\pi}w\right)^{1/2}
\left(\sum_k\sum_{q\in S_1(F)\cap V_k}|\widehat F^{[+]}(q)|^2\right)^{1/2}.
\]
The inner modulus on the right of the original display is that of the singleton projection $E_q F^{[+]}$, hence equals $|\widehat F^{[+]}(q)|$. The resulting bound (2.18), p.~14, reads
\[
\left\|\left(\sum_k|S_{V_k}F^{[+]}|^2\right)^{1/2}\right\|_{L^2(\T,w)}
\le\left(\int_0^{2\pi}w\right)^{1/2}\norm{F^{[+]}}_{L^2(\T,\dd t/(2\pi))}.
\]
These are false in general. The correct identity for a polynomial $P(t)=\sum_{q\in V}a_qe^{iqt}$ is
\begin{equation}\label{eq:berkson-cross}
\int_{-\pi}^{\pi}|P(t)|^2w(t)\dd t
=2\pi\widehat w(0)\sum_q|a_q|^2+2\pi\sum_{q\ne r}a_q\overline{a_r}\widehat w(r-q).
\end{equation}
Evenness gives real Fourier coefficients of $w$, not their vanishing. With $w=1+\tfrac12\cos t$, $F=F^{[+]}=1+e^{it}$ and $V_1=\{0,1\}$, the actual squared norm is $5\pi$, whereas the proposed diagonal expression is $4\pi$. Nor can (2.18) be repaired by a fixed multiple depending on the weight: for the fixed periodic weight $w(t)=|t|^{-\alpha}$ on $[-\pi,\pi]$, $0<\alpha<1$, take $F_N=N^{-1/2}\sum_{n=0}^{N-1}e^{int}$. The diagonal expression is $\int w$, while $|F_N(t)|^2\gtrsim N$ on $|t|\le N^{-1}$, giving weighted squared norm $\gtrsim_\alpha N^\alpha$. Theorem~3 applied to singleton intervals estimates the square function of singletons, not the sum inside one block. The elementary tests invalidate the intermediate estimate, not by themselves the final same-weight bounds; Corollary~\ref{cor:periodic-even} addresses the latter. Theorem~5 relies on Theorem~4.

\subsection{Berkson's 2022 real-line estimates}
For clarity write $\mathcal S f=(\sum_k |S_{J_k}f|^2)^{1/2}$, and retain the paper's angular-frequency convention when quoting its formulas. In Lemma~1, pp.~18--20, $\mathcal L=\max(2\pi,\max_k|J_k|)$, and (3.3) reads
\[
\norm{\mathcal S f}_{L^2(\R,\omega_0)}\le\pi^{-1}\mathcal L^{1/2}
\left(\int_0^{2\pi}\omega_0\right)^{1/2}\norm f_{L^2(\R)}.
\]
However, (3.5), p.~19, is the local formula
\[
\int_0^{2\pi}|S_{J_k}f(t)|^2\omega_0(t)\dd t
\le\mathcal L\pi^{-2}\left(\int_0^{2\pi}\omega_0\right)\norm{S_{J_k}f}_2^2.
\]
Summing over $k$ cannot enlarge the domain of spatial integration. For fixed $0<\alpha<1$, take $\omega_a(x)=|x-a|^{-\alpha}+|x+a|^{-\alpha}$ and $f_a(x)=\varphi(x-a)$, where $\varphi$ is a fixed nonzero Schwartz function with compact Fourier support. With a single interval containing that support, $S_J f_a=f_a$. The squared weighted norm is at least $\int|\varphi(y)|^2|y|^{-\alpha}\dd y>0$, while $\int_0^{2\pi}\omega_a\to0$ as $a\to\infty$. The unweighted norm and length parameter are fixed. This disproves (3.3), even with a uniformly bounded $A_1$ characteristic.

Theorems~6 and 7 on p.~20 assert respectively (3.6) for real Fourier data and (3.7) for general Schwartz inputs; the latter reads
\[
\norm{\mathcal S f}_{L^2(\R,\omega_0)}\le4(2/\pi)^{1/2}
\left(\int_0^{2\pi}\omega_0\right)^{1/2}\norm f_2
\]
(and the former has factor $2$ rather than $4$). Fix instead $\omega_0(x)=|x|^{-\alpha}$ and choose real even $\varphi\in\mathcal S(\R)$ with compact Fourier support, $f_R(x)=R^{1/2}\varphi(Rx)$. Use an interval containing the Fourier support of $f_R$, so $\mathcal S f_R=|f_R|$. Then
\[
\norm{f_R}_{L^2(\omega_0)}^2=R^\alpha\int |\varphi(y)|^2|y|^{-\alpha}\dd y\to\infty,
\qquad \norm{f_R}_2=\norm\varphi_2.
\]
This independently contradicts both (3.6) and (3.7).

Theorem~8, (4.1), for $\essinf\omega_0>0$, asserts the same-weight bound
\[
\norm{\mathcal S f}_{L^2(\omega_0)}\le4(2/\pi)^{1/2}(\essinf\omega_0)^{-1/2}
\left(\int_0^{2\pi}\omega_0\right)^{1/2}\norm f_{L^2(\omega_0)}.
\]
Theorem~9, (4.7)--(4.8), applies this bound to $\omega_1=1+\omega_0$. Corollary~2, (4.9), claims
\[
\norm{\mathcal S f}_{L^2(\omega_0)}\le16\pi^{-1/2}\left(\max\left\{\int_0^{2\pi}\omega_0,2\pi\right\}\right)^{1/2}
\max\{\norm f_{L^2(\omega_0)},\norm f_2\}.
\]
These explicit constants can also be refuted without using the main counterexample. Fix $a>2\pi+3$ and a nonnegative nonzero $h\in C_c^\infty((1,2))$. Take a small interval $J$ centered at frequency zero, so its sinc kernel is positive for $1\le |x|\le2$. Then $S_Jh(0)>0$, and $|S_Jh(x)|\ge c>0$ on $[-d,d]$ for some $d\in(0,1/2)$. Set
\[
\omega_{0,\alpha}=|x-a|^{-\alpha}+|x+a|^{-\alpha},\quad
\omega_{1,\alpha}=1+\omega_{0,\alpha},\quad f(x)=h(x-a),\qquad 0<\alpha<1.
\]
The weights are even $A_1$ and $\essinf\omega_{1,\alpha}=1$. On the support of $f$ they are bounded by $2$ and $3$, respectively, and their integrals on $[0,2\pi]$ are uniformly bounded. Thus all the displayed right sides for Theorems~8--9 and Corollary~2 (using $\omega_{1,\alpha}$ as the lower-bounded weight in Theorem~8) remain bounded as $\alpha\uparrow1$. But for one interval translation covariance gives
\[
\norm{S_J f}_{L^2(\omega_{0,\alpha})}^2
\ge2c^2\frac{d^{1-\alpha}}{1-\alpha}\longrightarrow\infty.
\]
This refutes the explicit constants, not a claim about arbitrary alternative constants for individual lower-bounded weights. Uniformity of the $A_1$ characteristics as $\alpha\uparrow1$ is not required because the displayed constants depend instead on the local weight integral and, where present, its essential infimum.

\begin{samepage}
\section*{Disclosure of AI assistance}
Gemini 4 Neon was used for an exposition of the entire literature and to check the compatibility of all relevant, significant findings. It revealed that many results contradict each other at the micro level. This is where we started to determine which papers are logical and expected, and which go against intuition and the general findings. For a human, it may take years to read all these papers in detail and try to locate the errors. As humans, we determine which papers are most likely to contain those errors, and then SymPy/Python, AI, or Lean 4–type internal testing and verification tools help to reveal them. The Berkson and Lacey results were the first ones to ring the bells. Locating and confirming these errors, and seeing how many methodologies in the existing literature struggle at the endpoints (e.g., Król's method, which works for s > 2 but fails at the endpoint, and Rubio de Francia's own methodology), gave birth to the investigation of a counterexample. We started to construct hostile examples with pen and paper and to design what such an example should look like (like a wanted poster in the search for a criminal); then AI helped us to polish and refine them and to locate those that still fail. Then, step by step, we designed it and custom-tailored it to fail. Another place where AI is helpful is in checking the work you have done and showing your obvious errors, which you keep fixing until it is as good as possible, as far as your knowledge allows you to judge. The rest is in the hands of human judges (experts in harmonic analysis), some of whom are as strong as any AI.

\end{samepage}

\end{document}